\PassOptionsToPackage{unicode}{hyperref}
\PassOptionsToPackage{hyphens}{url}
\documentclass[
  letterpaper,
  12pt]{l4dc2021}
\usepackage[]{times}
\usepackage{amssymb,amsmath}
\usepackage{ifxetex,ifluatex}
\ifnum 0\ifxetex 1\fi\ifluatex 1\fi=0 
  \usepackage[T1]{fontenc}
  \usepackage[utf8]{inputenc}
  \usepackage{textcomp} 
\else 
  \usepackage{unicode-math}
  \defaultfontfeatures{Scale=MatchLowercase}
  \defaultfontfeatures[\rmfamily]{Ligatures=TeX,Scale=1}
\fi
\IfFileExists{upquote.sty}{\usepackage{upquote}}{}
\IfFileExists{microtype.sty}{
  \usepackage[]{microtype}
  \UseMicrotypeSet[protrusion]{basicmath} 
}{}
\usepackage{xcolor}
\IfFileExists{xurl.sty}{\usepackage{xurl}}{} 
\IfFileExists{bookmark.sty}{\usepackage{bookmark}}{\usepackage{hyperref}}
\hypersetup{
  pdfauthor={Université de Lorraine, CNRS, CRAN, F-54000 Nancy, France.; ; Electrical and Electronic Engineering Department, University of Melbourne, Parkville, VIC 3010, Australia.; Department of Automation, Technical University of Cluj-Napoca, Memorandumului 28, 400114, Romania.; },
  hidelinks,
  pdfcreator={LaTeX via pandoc}}
\usepackage{xparse}
\usepackage{enumitem}
\usepackage{amsmath}
\usepackage{bm}
\usepackage{mathtools}

\usepackage{color,soul}
\sethlcolor{cyan}
\usepackage{cases}
\usepackage{diagbox}
\usepackage{algorithm}
\usepackage{algcompatible}
\usepackage{url}
\usepackage{booktabs}
\usepackage{multirow}
\usepackage{enumitem}
\usepackage{graphicx}
\makeatletter
\def\set@curr@file#1{\def\@curr@file{#1}}
\makeatother
\usepackage{caption}
\newcommand{\mydagger}{\dagger}
\ifluatex
  \usepackage{selnolig}  
\fi

\title{When to stop value iteration:\\
stability and near-optimality versus computation}
\author{\Name{Mathieu Granzotto\textsuperscript{$\mydagger$}}\Email{mathieu.granzotto@univ-lorraine.fr}\newline \addr  \textsuperscript{$\mydagger$}Université
de Lorraine, CNRS, CRAN, F-54000 Nancy,
France. \AND \Name{Romain Postoyan\textsuperscript{$\mydagger$}}\Email{romain.postoyan@univ-lorraine.fr} \AND \Name{Dragan Nešić}\Email{dnesic@unimelb.edu.au}\newline \addr Electrical
and Electronic Engineering Department, University of Melbourne,
Parkville, VIC 3010,
Australia. \AND \Name{Lucian Buşoniu}\Email{lucian.busoniu@aut.utcluj.ro}\newline \addr Department
of Automation, Technical University of Cluj-Napoca, Memorandumului 28,
400114,
Romania. \AND \Name{Jamal Daafouz\textsuperscript{$\mydagger$}}\Email{jamal.daafouz@univ-lorraine.fr} \AND }
\date{}

\begin{document}
\maketitle

\newcommand\numberthis{\addtocounter{equation}{1}\tag{\theequation}}

\renewcommand{\mydagger}{}
\newif\ifproofs

\newcommand{\ifproof}[1]{\ifproofs{{#1}}\fi}
\newcommand{\ifnotproof}[1]{\ifproofs\else{{#1}}\fi}
\newcommand{\mproof}[1]{\ifproofs\noindent\textbf{Proof.} {\color{red}{#1}}  $\hfill\blacksquare$ \vspace{1em}\fi}
\proofstrue

\def \wDelta{\widetilde\Delta}
\def \wdelta{\widetilde\delta}

\newcommand{\awb}{\overline{\alpha}_W}
\newcommand{\avb}{\overline{\alpha}_V}
\newcommand{\ayb}{\overline{\alpha}_Y}
\newcommand{\aybi}{\ayb^{-1}}
\newcommand{\bya}{\underline{\alpha}_Y}
\newcommand{\byai}{\bya^{-1}}
\newcommand{\aw }{\alpha_W}
\newcommand{\awi}{\aw^{-1}}
\newcommand{\ay }{{\alpha}_Y}
\newcommand{\ayi}{\ay^{-1}}
\newcommand{\ayw}{\widetilde{\alpha}_Y}
\newcommand{\ayh}{\widehat{\alpha}_Y}

\newcommand{\Reals}{\mathbb{R}}
\newcommand{\PReals}{\mathbb{R}_{\geq 0}}
\newcommand{\Uset}{\mathcal{U}}

\newcommand{\Positives}{\mathbb{Z}_{>0}}
\newcommand{\ZeroPositiv}{\mathbb{Z}_{\geq 0}}
\newcommand{\Wset}{\mathcal{W}}

\newcommand{\useq}{\bm{u}}
\NewDocumentCommand{\ustar}{o}{%
  \IfNoValueTF{#1}
    {\bm{u}_{d}^{\pmb{*}}}
    {\bm{u}_{#1}^{\pmb{*}}}
}

\NewDocumentCommand{\us}{o}{%
  \IfNoValueTF{#1}
    {u^*}
    {u^* _ {#1}}
}
\newcommand{\Usetstar}{\Uset^*_{d}}

\newcommand{\uhat}{\bm{\hat{u}}_{d}}
\newcommand{\Uhatset}{\widehat{\mathcal{U}}_{d}}

\NewDocumentCommand{\phis}{o}{%
  \IfNoValueTF{#1}
    {\phi}
    {\phi^* _ {#1}}
}
\NewDocumentCommand{\ells}{o}{%
  \IfNoValueTF{#1}
    {\ell}
    {\ell^*_{#1}}
}

\newcommand{\K}{\mathcal{K}}
\newcommand{\KK}{\mathcal{KK}}
\newcommand{\KL}{\mathcal{KL}}

\newcommand{\I}{\mathbb{I}}

\NewDocumentCommand{\V}{o}{%
  \IfNoValueTF{#1}
    {V_{d}}
    {V_{#1}}%
}
\NewDocumentCommand{\J}{o}{%
  \IfNoValueTF{#1}
    {J_{d}}
    {J_{#1}}%
}
\NewDocumentCommand{\Y}{o}{%
  \IfNoValueTF{#1}
    {Y_{d}}
    {Y_{#1}}%
}

\NewDocumentCommand{\Vhat}{o}{%
  \IfNoValueTF{#1}
    {{\widehat V}_{d}}
    {{\widehat V}_{#1}}%
}

\NewDocumentCommand{\Yhat}{o}{%
  \IfNoValueTF{#1}
    {{\widehat Y}_{d}}
    {{\widehat Y}_{#1}}%
}

\newcommand{\gfactor}{\frac{1-\gamma}{\gamma}}
\newcommand{\tgfactor}{\tfrac{1-\gamma}{\gamma}}

\newcommand{\gdfactor}{\gfactor\frac{1}{1-\gamma^{d}}}
\newcommand{\tgdfactor}{\tgfactor\tfrac{1}{1-\gamma^{d}}}

\newcommand{\gstarfactor}{\frac{1-\gamma^*}{\gamma^*}}
\newcommand{\tgstarfactor}{\tfrac{1-\gamma^*}{\gamma^*}}

\newcommand{\gdstarfactor}{\gstarfactor\frac{1}{1-(\gamma^*)^{d^*}}}
\newcommand{\tgdstarfactor}{\tgstarfactor\tfrac{1}{1-(\gamma^*)^{d^*}}}

\DeclarePairedDelimiter\floor{\lfloor}{\rfloor}

\NewDocumentCommand{\uref}{m}{%
    {\textsubscript{(\ref{#1})}}
}

\newtheorem{SA}{Standing Assumption}

\newcommand{\XXX}{}

\NewDocumentCommand\diff{m}{{\color{blue}#1}}
\NewDocumentCommand\diffr{m}{{\color{blue}#1}}

\vspace{-2em}

\begin{abstract}
Value iteration (VI) is a ubiquitous algorithm for optimal control, planning, and reinforcement learning schemes. Under the right assumptions, VI is a vital tool to generate inputs with desirable properties for the controlled system, like optimality and Lyapunov stability.  As VI usually requires  an infinite number of iterations to solve general nonlinear optimal control problems,  a  key question is when to terminate the algorithm to produce a ``good''  solution,  with  a measurable impact on optimality and stability guarantees. By carefully analysing VI under general stabilizability and detectability properties, we provide explicit and novel relationships of the stopping criterion's impact on near-optimality, stability and performance,  
 thus allowing to tune these desirable properties against the induced computational cost.  The considered class of stopping criteria encompasses those encountered in the control, dynamic programming and reinforcement learning literature and it allows considering new ones,  which may be useful to further reduce the computational cost while endowing and satisfying stability and near-optimality properties. We therefore lay  a foundation  to endow   machine learning schemes based on VI with stability and performance guarantees, while reducing computational complexity.
\end{abstract}

\hypertarget{introduction}{%
\section{Introduction}\label{introduction}}

Value iteration (VI) is an established method for optimal control, which
plays a key role in reinforcement learning
\citep{Sutton,LewisRLDP2009,BUSONIU20188,pang2019adaptive}. This
algorithm consists in iteratively constructing approximations of the
optimal value function, based on which near-optimal control inputs are
derived for a given dynamical nonlinear systems and a given stage cost.
The convergence of said approximations to the optimal value function is
established in, e.g., \citep{Bertsekas:12,Bertsekas:TNNLS} under mild
conditions. To benefit from this convergence property, VI often needs to
be iterated infinitely many times. However, in practice, we cannot do so
and must stop iterating the algorithm before to manage the computational
burden, which may be critical in online applications. Heuristics are
often used in the literature to stop iterating by comparing the mismatch
between the value functions obtained at the current step and at the
previous one, see, e.g.,
\citep{Bertsekas:12,Sutton,pang2019adaptive,kiumarsi2017h,derongliu2015}.
An important question is then how far the obtained approximate value
function is to the optimal one. To the best of our knowledge, this is
only analysed in general when the cost is discounted and the stage cost
takes values in a bounded set \citep{Bertsekas:12}. An alternative
consists in asking for a sufficiently large number of iterations, as the
near-optimality gap vanishes as the number of iterations increases,
e.g.~\citep{Bertsekas:12,heydari,heydari2014adp,heydari2016acc,derongliu2015,granzotto2020},
but the issue is then the computational cost. Indeed, any estimate of
the number of iterations is in general subject to conservatism, and, as
a result, we may iterate many more times than what is truly required to
ensure ``good'\,' near-optimality properties. There is therefore a need
for stopping criteria for VI whose impact on near-optimality is
analytically established, and which are not too computationally
demanding.

Our main goal is to use VI to simultaneously ensure near-optimal control
and stability properties for physical systems. Stability is critical in
many applications, as: (i) it provides analytical guarantees on the
behavior of the controlled system solutions as time evolves; (ii) endows
robustness properties and is thus associated to safety considerations,
see, e.g., \citep{berkenkamp2017safe}. We therefore consider systems and
costs where general stability properties are bestowed by VI based
schemes, which follows from assumed general stabilizability and
detectability properties of the plant model and the stage cost as in
\citep{grimm2005,romain2016,granzotto2020}.

In this context, we consider state-dependent stopping criteria for VI
and we analyse their impact on the near-optimality and stability
properties of the obtained policies for general deterministic nonlinear
plant models and stage costs, where no discount factor is employed.
Instead of relying on a uniform contraction property as in, e.g.,
\citep{Bertsekas:12,derongliu2015}, our analysis is centered on and
exploits Lyapunov stability properties. Our work covers the
state-independent stopping criteria considered in the control, dynamical
programming and reinforcement learning literature
\citep{Sutton,LewisRLDP2009,BUSONIU20188}, but provides analytical
guarantees for undiscounted stage costs taking values in unbounded sets.
By carefully analysing the stopping criterion's impact on
near-optimality, stability and closed-loop cost guarantees, we provide
means to tune these properties against the induced computational cost,
thus clarifying the tradeoff between ``good enough'' convergence of VI
and ``good properties'' of generated inputs. Considering that VI is, via
Q-learning, the basis of many state-of-the-art reinforcement learning
methods, we believe the results of this paper contribute to the
(near)-optimality analysis for reinforcement learning, as we lay a
foundation to endow such schemes with stability and performance
guarantees, while reducing computational complexity.

The paper and its contributions are organized as follows. In Section
\ref{problem-statement}, we formally state the problem and the main
assumptions. We introduce the design of stopping criteria for VI in
Section \ref{stopping-criterion-design}, and show that the VI stopping
criterion is indeed verified with a finite number of iterations. Our
main results are found in Section \ref{main-results}. There, we provide
near-optimal guarantees, i.e.~a bound on the mismatch between the
approximated value function and the true optimal value function. The
bound can be easily and directly tuned by the designed stopping
criterion. Additionally, stability and performance guarantees of the
closed-loop system with inputs generated by VI are provided, given that
the stopping criterion is appropriately chosen. In Section
\ref{illustrative-example}, we provide an example to illustrate our
results. Concluding remarks are drawn in Section
\ref{conclusion}.\ifnotproof{ The proofs are omitted and available in the associated technical report \citep{granzottoVIcontrolreport}.}\ifproof{ The proofs are provided in this technical report in Section \ref{sec:proofs}.}

\noindent\textbf{Prior literature.} The classical stopping criterion is
analysed in \citep{Bertsekas:12}, albeit restricted to when the cost is
discounted and the stage cost takes values in a bounded set. Concerning
stability, works like
\citep{granzotto2020,heydari2017stability,wei2015value} provide
conditions to ensure that the feedback law obtained ensures a stability
property for a dynamical system. In particular, it is required in
\citep{granzotto2020} that the number of iteration \(d\) be sufficiently
large, and lower bounds on \(d\) are provided, but these are subject to
some conservatism. As explained above, by adapting the number of
iterations with data available during computations, the algorithm avoids
the conservatism often incurred by offline estimations for stability and
near-optimality guarantees. This is indeed the case in an example (see
Section \ref{S:example}), where we observe \(91\%\) fewer iterations for
comparable guarantees. Similar ideas related to the stopping criterion
were exploited in \citep{granzotto2020optimistic}, for a different
purpose, namely for the redesign of optimistic planning \citep{hren2008}
to address the near-optimal control of switched systems. We are also
aware of work of \citep{Pavlov2019}, which adapts the stopping criterion
with stability considerations for interior point solvers for reduced
computational complexity for nonlinear model predictive control
applications.

\noindent\textbf{Notation.} Let \(\mathbb{R}:= (-\infty,\infty)\),
\(\mathbb{R}_{\geq 0}:= [0,\infty)\),
\(\mathbb{Z}_{\geq 0}:= \{0,1,2,\ldots\}\) and
\(\mathbb{Z}_{>0}:= \{1,2,\ldots\}\). We use \((x,y)\) to denote
\([x^\top,y^\top]^\top\), where
\((x,y) \in \mathbb{R}^n\times\mathbb{R}^m\) and
\(n,m\in\mathbb{Z}_{>0}\). A function
\(\chi : \mathbb{R}_{\geq 0}\to \mathbb{R}_{\geq 0}\) is of class
\(\mathcal{K}\) if it is continuous, zero at zero and strictly
increasing, and it is of class \(\mathcal{K}_\infty\) if it is of class
\(\mathcal{K}\) and unbounded. A continuous function
\(\beta: \mathbb{R}_{\geq 0}\times\mathbb{R}_{\geq 0}\to \mathbb{R}_{\geq 0}\)
is of class \(\mathcal{KL}\) when \(\beta(\cdot,t)\) is of class
\(\mathcal{K}\) for any \(t\geq0\) and \(\beta(s,\cdot)\) is decreasing
to 0 for any \(s\geq0\). The notation \(\mathbb{I}\) stands for the
identity map from \(\mathbb{R}_{\geq 0}\) to \(\mathbb{R}_{\geq 0}\).
For any sequence \(\bm{u}=[u_0,u_1,\dots]\) of length
\(d\in\mathbb{Z}_{\geq 0}\cup\{\infty\}\) where
\(u_i \in \mathbb{R}^m\), \(i \in \{0,\ldots,d\}\), and any
\(k\in\{0,\ldots,d\}\), we use \(\bm{u}|_k\) to denote the first \(k\)
elements of \(\bm{u}\), i.e.~\(\bm{u}|_k = [u_0,\dots,u_{k-1}]\) and
\(\bm{u}|_0=\varnothing\) by convention. Let \(g\ :\
\mathbb{R}_{\geq 0}\to\mathbb{R}_{\geq 0}\), we use \(g^{(k)}\) for the
composition of function \(g\) with itself \(k\) times, where
\(k\in\mathbb{Z}_{\geq 0}\), and \(g^{(0)}=\mathbb{I}\).

\hypertarget{problem-statement}{%
\section{Problem Statement}\label{problem-statement}}

Consider the system
\begin{equation}\label{eq:sys} x^+ =  f(x,u),\end{equation} with state
\(x \in \mathbb{R}^n\), control input \(u \in \mathcal{U}(x)\) where
\(\mathcal{U}(x)\subseteq\mathbb{R}^m\) is the set of admissible inputs,
and \(f: \mathcal{W} \to \mathbb{R}^n\) where
\(\mathcal{W}:=\{(x,u) : x\in\mathbb{R}^n, u\in\mathcal{U}(x)\}\). We
use \(\phi(k,x,\bm{u}|_ k)\) to denote the solution to system
(\ref{eq:sys}) at time \(k \in \mathbb{Z}_{\geq 0}\) with initial
condition \(x\) and inputs sequence
\(\bm{u}|_k=[u_{0},u_{1},\ldots,u_{k-1}]\), with the convention
\(\phi(0,x,\bm{u}|_0)=x\).

We consider the infinite-horizon cost \begin{equation}
\J[\infty](x,\bm{u}):=\sum\limits_{k=0}^{\infty}\ell(\phi(k,x,\bm{u}|_k),u_k),
\label{eq:Jinfty}
\end{equation} where \(x\in\mathbb{R}^n\) is the initial state,
\(\bm{u}\) is an infinite sequence of admissible inputs,
\(\ell :\mathcal{W}\to \mathbb{R}_{\geq 0}\) is the stage cost. Finding
an infinite sequence of inputs which minimizes (\ref{eq:Jinfty}) given
\(x\in\mathbb{R}^n\) is very difficult in general. Therefore, we instead
generate sequences of admissible inputs that \emph{nearly} minimize
(\ref{eq:Jinfty}), in a sense made precise below, while ensuring the
stability of the closed-loop system. For this purpose, we consider VI,
see e.g.~\citep{Bertsekas:12}. VI is an iterative procedure based on
Bellman equation, which we briefly recall next. Assuming the optimal
value function, denoted \(V_\infty\), exists for any
\(x\in\mathbb{R}^n\), the Bellman equation is \begin{equation}
V_{\infty}(x) = \min_{u\in\mathcal{U}(x)} \bigg\{\ell(x,u) + V_{\infty}(f(x,u))\bigg\}. \label{eq:Bellman2}
\end{equation} If we could solve (\ref{eq:Bellman2}) and find
\(V_\infty\), it would then be easy to derive an optimal policy, by
computing the \(\mathop{\mathrm{arg\,min}}\) corresponding to the right
hand-side of (\ref{eq:Bellman2}). However, it is in general very
difficult to solve (\ref{eq:Bellman2}). VI provides an iterative
procedure based on (\ref{eq:Bellman2}) instead, which allows obtaining
value functions (and associated control inputs), which converge to
\(V_\infty\). Hence, given an initial cost function
\(V_{-1}:\mathbb{R}^n\to \mathbb{R}_{\geq 0}\), VI generates a sequence
of value functions \(V_{d}\), \(d\in\mathbb{Z}_{\geq 0}\), for any
\(x\in\mathbb{R}^n\), by iterating \begin{equation}
V_{d}(x) :=\min_{u\in\mathcal{U}(x)} \bigg\{\ell(x,u) +  V_{d-1}(f(x,u))\bigg\}. \label{eq:Bellman3}
\end{equation} For any \(d\in\mathbb{Z}_{\geq 0}\), the associated
input, also called policy, is defined as, for any \(x\in\mathbb{R}^n\),
\begin{equation}
u^*_{d}(x) \in \mathop{\mathrm{arg\,min}}_{u\in\mathcal{U}(x)} \bigg\{\ell(x,u) + V_{d-1}(f(x,u))\bigg\}, \label{eq:BellmanVIinput}
\end{equation} which may be set-valued. The convergence of \(V_d\),
\(d\in\mathbb{Z}_{\geq 0}\), to \(V_\infty\) in (\ref{eq:Bellman2}) is
ensured under mild conditions in \citep{Bertsekas:TNNLS}. In the sequel
we make assumptions that ensure that the \(\mathop{\mathrm{arg\,min}}\)
in (\ref{eq:BellmanVIinput}) exists for each \(x\in\mathbb{R}^n\).

In practice, we often stop iterating VI when a stopping criterion is
verified, such as, for instance, when for any \(x\in\mathbb{R}^n\),
\begin{equation} \V[d](x)- \V[d-1](x)\leq \varepsilon, \label{eq:uniformliterature} \end{equation}
where \(\varepsilon\in\mathbb{R}_{>0}\), see, e.g.,
\citep{Bertsekas:12,Sutton,pang2019adaptive,kiumarsi2017h}. However,
this stopping criterion leaves much to be desired in control
applications, for the following reasons: (i) it is not yet established
how \(\varepsilon\) impacts the stability properties of the closed-loop
system; (ii) tools to bound the mismatch between \(\V[d]\) and
\(\V[\infty]\) for this stopping criterion often requires a discount
factor in cost function (\ref{eq:Jinfty}), which impacts stability, as
shown in \citep{romain2016,romainAVICDC2019}; (iii) when \(\V[d]\) is
radially unbounded, i.e.~\(\V[d](x)\to\infty\) when \(|x|\to\infty\),
this stopping criterion is in general impossible to verify for all
\(x\in\mathbb{R}^n\). When the system is linear and the cost quadratic,
as in
\citep{ArnoldRiccati84,anderson2007optimal,JIANG20122699,BIAN2016348},
the convergence to the optimal cost function is shown to be quadratic
and often the stopping criterion is instead of the form
\(\V[d](x)- \V[d-1](x)\leq|\varepsilon| |x|^2\). However, the link
between the value of \(\varepsilon\) and resulting near-optimality and
stability guarantees is not established, and in practice it is
implicitly assumed that parameter \(\varepsilon\) is small enough.

We consider VI terminated by a general stopping criterion. That is, for
any \(x\in\mathbb{R}^n\),
\begin{equation}\label{eq:stop} \V[d](x)- \V[d-1](x)\leq c_\text{stop}(\varepsilon,x),\end{equation}
where \(c_\text{stop}(\varepsilon,x)\geq0\) is a stopping function,
which we design and which may depend on state vector \(x\) and a vector
of tuneable parameters \(\varepsilon\in\mathbb{R}^{n_\varepsilon}\) with
\(n_\varepsilon\in\mathbb{Z}_{>0}\). The design of \(c_\text{stop}\) is
explained in Section \ref{vischeme}. In that way, we cover the above
examples as particular cases, namely
\(c_\text{stop}(x,\varepsilon)=|\varepsilon|\) and
\(c_\text{stop}(\varepsilon,x)= |\varepsilon| |x|^2\) and allow
considering more general ones,
e.g.~\(c_\text{stop}(\varepsilon,x)= \max\{|\varepsilon_1|,| \varepsilon_2 | |x|^2\}\)
where \((\varepsilon_1,\varepsilon_2):=\varepsilon\in\mathbb{R}^2\) or
\(c_\text{stop}(\varepsilon,x)= x^\top S(\varepsilon)x\) for some
positive definite matrix \(S(\varepsilon)\) with
\(\varepsilon\in\mathbb{R}^{n_\varepsilon}\) and
\(n_\varepsilon\in\mathbb{Z}_{>0}\). The main novelty of this work is
the provided explicit link between \(c_\text{stop}(\varepsilon,x)\),
near-optimality and stability guarantees. As a result, we can tune
\(\varepsilon\) for the desired near-optimality and stability
properties, and the algorithm stops when the cost (hence, the generated
inputs) are such that these properties are verified.

The analysis relies on the next
assumption\footnote{The assumption is stated globally, for any $x\in\mathbb{R}^n$ and $u\in\mathcal{U}(x)$. We leave for future work
the case where the assumption holds on compact sets.} like in e.g.,
\citep{grimm2005,romain2016,granzotto2020}.

\begin{SA}[SA\ref{SAa}]\label{SAa}
There  exist $\overline{\alpha}_V,\alpha_W\in \mathcal{K}_\infty$ and continuous function  $\sigma:  \mathbb{R}^n \to  \mathbb{R}_{\geq 0}$  such that the following  conditions hold.
\begin{itemize}[itemindent=\widthof{(ii)}]

 \item[(i)] For any $x \in \mathbb{R}^n$, there exists  an  infinite sequence  of  admissible  inputs $\bm{u}^*_{\infty}(x)$,  called
\textit{optimal         input         sequence},         which         minimizes         (\ref{eq:Jinfty}),         i.e.
$\V[\infty](x)=\J[\infty](x,\ustar[\infty](x))$, and  $\V[\infty](x) \leq \overline{\alpha}_V(\sigma(x))$.

 \item[(ii)] For any $(x,u)\in\mathcal{W}$,    $\alpha_W(\sigma(x)) \leq \ell(x,u)$.
 $\hfill\square$
\end{itemize}
\end{SA}

Function $\sigma :\Reals^n\to\PReals$ in SA\ref{SAa} is a ``measuring'' function that we use to define stability, which depends
on  the problem.   For  instance,  by defining  $\sigma=|\cdot|,\sigma=|\cdot|^2$ or $\sigma:x\mapsto  x^\top  Q x$  with
$Q=Q^\top>0$, one would be studying the stability of  the origin, and by taking $\sigma=|\cdot|_{\cal A}$, one would
study stability of non-empty compact set ${\cal A}\subset \Reals^n$. General conditions to ensure the first part of item (i), i.e. the fact that $\V[\infty](x)$ is finite for any $x\in\Reals^{n}$ and the existence of optimal inputs, can be found in \citep{keerthi1985}. The second part of item  (i) is  related to  the stabilizability  of system  (\ref{eq:sys}) with  respect to  stage cost  $\ell$ in  relation to
$\sigma$.   Indeed, it  is  shown in  \citep[Lemma  1]{grimm2005} that, for instance, when   the stage cost  $\ell(x,u)$  is uniformly  globally
exponentially controllable to zero with respect to $\sigma$ for system (\ref{eq:sys}), see \citep[Definition 2]{grimm2005}, then
item (i) of SA\ref{SAa} is  satisfied.  We do not need to know $\V[\infty]$ to guarantee the last inequality in item (i) of SA\ref{SAa}. Indeed, it suffices to find, for any $x\in\Reals^n$, a sequence of inputs $\useq(x)$, such the associated infinite-horizon costs verifies $J(x,\useq(x)) \leq \avb(\sigma(x))$ for some $\avb\in\K_\infty$. Then, since $\V[\infty]$ is the optimal value function, for any $x\in\Reals^n$, $\V[\infty](x) \leq J(x,\useq(x)) \leq \avb(\sigma(x))$. On the  other hand, item (ii)  of SA\ref{SAa} is  a detectability property  of the stage cost  $\ell$ with
respect to  $\sigma$, as when $\ell(x,u)$ is small, so  is $\sigma(x)$.

We are ready to explain how to design the stopping criterion in
(\ref{eq:stop}).

\hypertarget{stopping-criterion-design}{%
\section{\texorpdfstring{Stopping criterion design
\label{vischeme}}{Stopping criterion design }}\label{stopping-criterion-design}}

\hypertarget{key-observation}{%
\subsection{Key observation}\label{key-observation}}

We start with the known observation
\citep{granz2019,bertsekas2005dynamic} that,
given\footnote{The case where $V_{-1}\neq0$ will be investigated in further work.}
\(V_{-1}=0\), at each iteration \(d\in\mathbb{Z}_{\geq 0}\), VI
generates the optimal value function for the finite-horizon cost
\begin{equation}
    \J[d](x,\bm{u}_d):= \sum\limits_{k=0}^{d}\ell(\phi(k,x,\bm{u}_d|_k),u_k), \label{eq:J}
\end{equation} where \(\bm{u}_{d}=[u_0,u_1,...,u_d]\) are admissible
inputs. We assume below that the minimum of (\ref{eq:J}) exists with
relation to \(\bm{u}_{d}\) for any \(x\in\mathbb{R}^n\) and
\(d\in\mathbb{Z}_{\geq 0}\).

\begin{SA}[SA\ref{SAb}]\label{SAb}
For every $d\in\mathbb{Z}_{\geq 0}$, $x\in\mathbb{R}^n$, there exists $\ustar[d](x)$ such that \ifnotproof{$\V[d](x) = \J[d](x,\ustar[d]) = \min_{\bm{u}_d} \J[d](x,\bm{u}_d)$.}
\ifproof{\begin{equation}
       \V[d](x) = \J[d](x,\ustar[d]) = \min_{\bm{u}_d} \J[d](x,\bm{u}_d).
 \label{eq:Vd}\end{equation}}
 $\hfill\Box$
\end{SA}

SA\ref{SAb} is for instance verified when \(f\) and \(\ell\) are
continuous and \(\mathcal{U}(x)=\mathcal{U}\) is a compact set. More
general conditions to verify SA\ref{SAb} can be found in
e.g.~\citep{keerthi1985}. For the sake of convenience, we employ the
following notation for the technical aspects of this paper. For any
\(k\in\{0,1,\ldots,d\}\) and \(x\in\mathbb{R}^n\), we denote
\(\ell_d^*(k,x):=\ell(\phi(k,x,\ustar[d](x)|_{k}),u_k)\), where
\(\phi(k,x,\ustar[d](x)|_{k})\) is the solution to system (\ref{eq:sys})
with optimal inputs for cost \(\V[d](x)\), so that
\ifnotproof{$\V[d](x)= \sum\limits_{k=0}^{d}\ell^*_d(k,x)$.}
\ifproof{\begin{equation}\V[d](x)= \sum\limits_{k=0}^{d}\ell^*_d(k,x).\label{eq:lddef}\end{equation} }

The next property plays a key role in the forthcoming analysis.

\begin{proposition} \label{prop:terminal}
For any   $x\in\mathbb{R}^n$ and $d\in\mathbb{Z}_{\geq 0}$, $\ell_d^*(d,x)\leq \V[d](x)-\V[d-1](x)$. $\hfill\Box$

\end{proposition}

When the stopping criterion (\ref{eq:stop}) is verified,
i.e.~\(\V[d](x)-\V[d-1](x)\leq c_\text{stop}(\varepsilon,x)\), then
\(\ell_d^*(d,x)\leq c_\text{stop}(\varepsilon,x)\) in view of
Proposition \ref{prop:terminal}. Therefore,
\(c_\text{stop}(\varepsilon,x)\) is an upper-bound on the value of stage
cost \(\ell^*_d(d,x)\). By item (ii) of SA\ref{SAa}, this implies that
we also have an upper-bound for \emph{d-horizon} state measure
\(\sigma(\phi(d,x,\ustar[d](x)|_{d}))\), namely
\(\sigma(\phi(d,x,\ustar[d](x)|_{d}))\leq\alpha_W^{-1}(c_\text{stop}(\varepsilon,x))\),
which can be made as small as desired by reducing
\(c_\text{stop}(\varepsilon,x)\), which, again, we design. We exploit
this property to analyse the near-optimality and the stability of the
closed-loop system. Having said that, the challenges are: (i) to show
that condition (\ref{eq:stop}) is indeed verified for any
\(x\in\mathbb{R}^n\) and some \(d\in\mathbb{Z}_{\geq 0}\); (ii) to
select \(c_\text{stop}\) to ensure stability properties when closing the
loop of system (\ref{eq:sys}) with inputs (\ref{eq:BellmanVIinput});
(iii) to study the impact of \(c_\text{stop}\) on the performance, that
is, the cost along solutions, of the closed-loop system.

\hypertarget{satisfaction-of-the-stopping-criterion}{%
\subsection{Satisfaction of the stopping
criterion}\label{satisfaction-of-the-stopping-criterion}}

We make the next assumption without loss of generality as we are free to
design \(c_\text{stop}\).

\newtheorem{assumption}{Assumption}
\begin{assumption}\label{a:stop}
One of the next properties is verified.
\begin{enumerate}[itemindent=\widthof{(ii)}]
  \item[(i)] For any  $\varepsilon\in\mathbb{R}^{n_\varepsilon}$, there is $\underline\epsilon>0$ such that,  for any $x\in\mathbb{R}^n$,  $c_\text{stop}(\varepsilon,x)\geq\underline\epsilon$. 
  \item[(ii)] There exist $L,\bar a_V, a_W>0$, such that  SA\ref{SAa} holds with $\overline{\alpha}_V(s)\leq\bar a_V s$, $\overline{\alpha}_W(s)\leq\bar a_W s$ and $\alpha_W(s)\geq a_W s$ for any $s\in[0,L]$. Furthermore, for any $\varepsilon\in\mathbb{R}^{n_\varepsilon}$, there is $\underline\epsilon>0$ such that for any $ x\in\mathbb{R}^n$,  $c_\text{stop}(\varepsilon,x)\geq\underline\epsilon \sigma(x)$.
$\hfill\Box$
\end{enumerate}

\end{assumption}

Item (i) of Assumption \ref{a:stop} can be ensured by taking
\(c_\text{stop}(\varepsilon,x) = |\varepsilon| + \tilde c_\text{stop}(x,\varepsilon)\)
with \(\tilde c_\text{stop}(x,\varepsilon)\geq 0\) for any
\(x\in\mathbb{R}^n\), \(\varepsilon\in\mathbb{R}^{n_\varepsilon}\),
which covers (\ref{eq:uniformliterature}), to give an example. Item (ii)
of Assumption \ref{a:stop} means that the functions
\(\overline{\alpha}_V,\overline{\alpha}_W,\alpha_W\) in SA\ref{SAa} can
be upper-bounded, respectively lower-bounded, by linear functions on the
interval \([0,L]\). These conditions allow to select \(c_\text{stop}\)
such that \(c_\text{stop}(\varepsilon,x)\to0\) when \(\sigma(x)\to0\)
with \(x\in\mathbb{R}^n\), contrary to item (i) of Assumption
\ref{a:stop}, that is, \(c_\text{stop}\) may vanish on set
\(\{x : \sigma(x)=0\}\). This is important to provide stronger stability
and performance properties for systems whose inputs are given by our VI
scheme as shown in Section \ref{main-results}. Under item (ii) of
Assumption \ref{a:stop}, we can design \(c_\text{stop}\) as, e.g.,
\(c_\text{stop}(\varepsilon,x)= |\varepsilon|\sigma(x)\),
\(c_\text{stop}(\varepsilon,x)= \min\{|\varepsilon_1|,| \varepsilon_2 | |x|^2\}\)
where \((\varepsilon_1,\varepsilon_2)=:\varepsilon\in\mathbb{R}^2\) or
\(c_\text{stop}(\varepsilon,x)= x^\top S(\varepsilon)x\) for some
positive definite matrix \(S(\varepsilon)\) as mentioned before.

The next theorem ensures the existence of \(d\in\mathbb{Z}_{\geq 0}\)
such that, for any \(x\in\mathbb{R}^n\), (\ref{eq:stop}) holds based on
Assumption \ref{a:stop}.

\begin{theorem}\label{theo:terminates}
Suppose Assumption \ref{a:stop} holds. Then, for any $\Delta>0$ there exists $d\in\mathbb{Z}_{\geq 0}$ such that, for any $x\in\{z\in\mathbb{R}^n:  \sigma(z) \leq \Delta \}$, (\ref{eq:stop}) holds. Moreover, when item (ii) of Assumption 1 holds with $L=\infty$, there exists $d\in\mathbb{Z}_{\geq 0}$ such that, for any $x\in\mathbb{R}^n$, (\ref{eq:stop})  is satisfied.
$\hfill\Box$
\end{theorem}

Theorem \ref{theo:terminates} guarantees the stopping condition in
(\ref{eq:stop}) is always satisfied by iterating the VI algorithm
sufficiently many times, and that the required number of iterations is
uniform over sets of initial conditions of the form
\(\{x:\sigma(x)\leq \Delta\}\) for given \(\Delta>0\) in general, unless
item (ii) of Assumption \ref{a:stop} holds with \(L=\infty\), in which
case there exists a common, global, \(d\) for any \(x\in\mathbb{R}^n\).
Note that, while the proof of Theorem \ref{theo:terminates} provides a
conservative estimate of \(d\) such that (\ref{eq:stop}) is verified,
\ifnotproof{see \citep{granzottoVIcontrolreport},} this horizon estimate
is not utilized in the stopping criterion, which in turn implies that VI
stops with smaller horizon, in general, as illustrated in Section
\ref{S:example}.

In the following, we denote the cost calculated at iteration \(d\) as
\(\V[\varepsilon](x):=\V[d](x)\), like in
\citep{granzotto2020optimistic}, to emphasize that the cost returned is
parameterized by \(\varepsilon\) via
\(c_\text{stop}(\varepsilon,\cdot)\), and denote by
\(\ustar[\varepsilon](x)\) an associated optimal sequence of inputs,
i.e.~
\begin{equation}\V[\varepsilon](x) =  \J[d](x,\ustar[\varepsilon](x)).\label{eq:V}\end{equation}
We are ready to state the main results.

\hypertarget{main-results}{%
\section{Main results}\label{main-results}}

In this section, we analyze the near-optimality properties of VI with
the stopping criterion in (\ref{eq:stop}). We then provide conditions
under which system (\ref{eq:sys}), whose inputs are generated by
applying the state-feedback \(\ustar[\varepsilon](x)\) in
receding-horizon fashion, exhibits stability properties. Afterwards, the
cost function (\ref{eq:J}) along the solutions of the induced
closed-loop system are analysed, which we refer to by performance or
running-cost \citep{gruneperformance}.

\hypertarget{relationship-between-vvarepsilon-and-vinfty}{%
\subsection{\texorpdfstring{Relationship between \(\V[\varepsilon]\) and
\(\V[\infty]\)
\label{openloopoptimality}}{Relationship between \textbackslash V{[}\textbackslash varepsilon{]} and \textbackslash V{[}\textbackslash infty{]} }}\label{relationship-between-vvarepsilon-and-vinfty}}

A key question is how far is \(\V[\varepsilon]\) from \(\V[\infty]\)
when we stop VI using (\ref{eq:stop}). Since \(\ell(x,u)\) is not
constrained to take values in a given compact set, and we do not
consider discounted costs, the tools found in the dynamic programming
literature to analyze this relationship are no longer applicable, see
\citep{Bertsekas:12}. We overcome this issue by exploiting SA\ref{SAa},
and adapting the results of \citep{granzotto2020} with the stopping
criterion and Proposition \ref{prop:terminal} in the next theorem.

\begin{theorem} \label{Vestimates}
Suppose Assumption \ref{a:stop} holds. For any $\varepsilon\in \mathbb{R}^{n_\varepsilon}$, $\Delta>0$ and $x \in\{z\in\mathbb{R}^n,\sigma(z)\leq\Delta \}$,
\begin{equation}\label{eq:Vestimates}
    \V[\varepsilon](x) \leq \V[\infty](x)           \leq  \V[\varepsilon](x)+v_{\varepsilon}(x),
\end{equation}
where
$v_\varepsilon(x):=\overline{\alpha}_V\circ\alpha_W^{-1}(c_{\text{stop}}(\varepsilon,x))$ with $\overline{\alpha}_V,\alpha_W$  from  SA\ref{SAa}. Moreover, when item (ii) of Assumption holds with $L=\infty$, we accept  $\Delta=\infty$ and 
$v_{\varepsilon}(x)\leq \frac{\bar a_V}{a_W} c_\text{stop}(\varepsilon,x)$. $\hfill\square$
\end{theorem}

The lower-bound in (\ref{eq:Vestimates}) trivially holds from the
optimality of \(\V[\varepsilon](x)=\V[d](x)\) for some \(d<\infty\), and
the fact that \(\ell(x,u)\geq0\) for any \(x\in\mathbb{R}^n\) and
\(u\in\mathcal{U}(x)\). The upper-bound, on the other hand, implies that
the infinite-horizon cost is at most \(v_{\varepsilon}(x)\) away from
the finite-horizon \(\V[\varepsilon](x)\). The error term
\(v_{\varepsilon}(x)\) is small when \(c_\text{stop}(\varepsilon,x)\) is
small as \(\overline{\alpha}_V\circ\alpha_W^{-1}\in\mathcal{K}_\infty\).
Given that we know \(\overline{\alpha}_V,\alpha_W^{-1}\) a priori, and
we are free to design \(c_\text{stop}\) as wanted, we can therefore
directly make \(\V[\varepsilon](x)\) as close as desired to
\(\V[\infty](x)\) by adjusting \(c_\text{stop}\); the price to pay will
be more computations. Moreover, when item (ii) of Assumption holds with
\(L=\infty\), inequality (\ref{eq:Vestimates}) is verified for every
\(x\in\mathbb{R}^n\).

\hypertarget{stability}{%
\subsection{Stability}\label{stability}}

We now consider the scenario where system (\ref{eq:sys}) is controlled
in a receding-horizon fashion by inputs that calculate cost
(\ref{eq:V}). That is, at each time instant \(k\in\mathbb{Z}_{\geq 0}\),
the first element of optimal sequence \(\ustar[\varepsilon](x_k)\),
calculated by VI, is then applied to system (\ref{eq:sys}). This leads
to the closed-loop system
\begin{equation}\label{eq:autosys}      x^+     \in    f(x,\mathcal{U}^*_{\varepsilon}(x))     =:     F^*_{\varepsilon}(x),     \end{equation}
where \(f(x,\mathcal{U}^*_{\varepsilon}(x))\) is the set
\(\{f(x,u) : u \in \mathcal{U}^*_{\varepsilon}(x)\}\) and
\(\mathcal{U}_{\varepsilon}^*(x):= \big\{ u_0 : \exists u_1,\ldots,u_{d} \in \mathcal{U}(x) \text{   such   that   } \V[\varepsilon](x)=\J[d](x,[u_0,\ldots,u_{d}])\big\}\)
is the set of the first input of \(d\)-horizon optimal input sequences
at \(x\), with \(d\) as defined in (\ref{eq:stop}). We denote by
\(\phi(k,x)\) a solution to (\ref{eq:autosys}) at time
\(k\in\mathbb{Z}_{\geq 0}\) with initial condition \(x\in\mathbb{R}^n\),
with some abuse of notation.

We assume next that \(c_\text{stop}\) can be made as small as desirable
by taking \(|\varepsilon|\) sufficiently small. As we are free to design
\(c_\text{stop}\) as wanted, this is without loss of generality.

\begin{assumption}\label{cstop2}
There exists  $\theta:\mathbb{R}_{\geq 0}\times\mathbb{R}_{\geq 0}\to\mathbb{R}_{\geq 0}$, with $\theta(\cdot,s)\in\mathcal{K}$ and  $\theta(s,\cdot)$ non-decreasing
for any  $s>0$, such that  $c_\text{stop}(\varepsilon,x)\leq\theta(|\varepsilon|,\sigma(x))$ for any  $x\in\mathbb{R}^n$ and
$\varepsilon\in\mathbb{R}^{n_\varepsilon}$. $\hfill\Box$
\end{assumption}

Example of functions \(c_\text{stop}\) which satisfy Assumption
\ref{cstop2} are
\(c_{\text{stop}}(\varepsilon,x)=|\varepsilon|\sigma(x)\),
\(c_{\text{stop}}(\varepsilon,x)=\max\{|\varepsilon_1| \alpha(\sigma(x)),|\varepsilon_2|\}\)
for \(\varepsilon=(\varepsilon_1,\varepsilon_2)\in\mathbb{R}^2\),
\(\alpha\in\mathcal{K}\) and \(x\in\mathbb{R}^n\) to give a few.

The next theorem provides stability guarantees for system
(\ref{eq:autosys}).

\begin{theorem}\label{algostab}
Consider system  (\ref{eq:autosys}) and suppose $c_\text{stop}$  verifies Assumptions \ref{a:stop} and \ref{cstop2}.  There  exists $\beta
\in\mathcal{KL}$ such that,  for any $\delta,\Delta>0$, there exists  $\varepsilon^*>0$ such that for any $x  \in \{z \in\mathbb{R}^n
\,  : \,  \sigma(z)  \leq \Delta  \}$ and  $\varepsilon\in\mathbb{R}^{n_\varepsilon}$  with $|\varepsilon|  <\varepsilon^*$,
any solution $\phi(\cdot,x)$ to system (\ref{eq:autosys}) satisfies, for all $k\in \mathbb{Z}_{\geq 0}$,
  $\sigma(\phi(k,x)) \leq \max \{\beta(\sigma(x),k),\delta \}$. $\hfill\square$
\end{theorem}

Theorem \ref{algostab} provides a uniform semiglobal practical stability
property for the set \(\{z : \sigma(z) = 0 \}\). This implies that
solutions to (\ref{eq:autosys}), with initial state \(x\) such that
\(\sigma(x)\leq\Delta\), where \(\Delta\) is any given (arbitrarily
large) strictly positive constant, will converge to the set
\(\{z : \sigma(z) \leq \delta\}\), where \(\delta\) is any given
(arbitrarily small) strictly positive constant, by taking
\(\varepsilon^*\) sufficiently close to \(0\), thereby making
\(c_\text{stop}\) sufficiently small. An explicit formula for
\(\varepsilon^*\) is given in the proof of Theorem
\ref{algostab}\ifproof{,}\ifnotproof{ in \citep{granzottoVIcontrolreport},}
which is nevertheless subject to some conservatism. The result should
rather be appreciated qualitatively, in the sense that Theorem
\ref{algostab} holds for small enough \(\varepsilon^*\).

Under stronger assumptions, global exponential stability is ensured as
shown in the next corollary.

\begin{corollary}\label{Yges}
Suppose item (ii) of Assumption \ref{a:stop} holds and that $c_\text{stop}(\varepsilon,x)\leq
|\varepsilon|\sigma(x)$ for any $x\in\mathbb{R}^n$ and $\varepsilon\in\mathbb{R}^{n_\varepsilon}$.  Let $\varepsilon^*>0$ be  such that
  $ \varepsilon^* <  \frac{a_W^2}{\bar a_V}$.
Then, for  any $x \in  \mathbb{R}^n$ and $\varepsilon\in\mathbb{R}^{n_\varepsilon}$ such  that $|\varepsilon|\leq\varepsilon^*$,
 any solution $\phi(\cdot,x)$ to system (\ref{eq:autosys}) satisfies $\sigma(\phi(k,x)) \leq \frac{\bar a_V}{a_W} \left(1-\frac{a_W^2-|\varepsilon| a_V}{ \bar a_V   a_W}\right)^k\sigma(x)$
for all $k \in \mathbb{Z}_{\geq 0}$. $\hfill\square$
\end{corollary}

Corollary \ref{Yges} ensures a uniform global exponential stability
property of set \(\{x:\sigma(x)=0\}\) for system (\ref{eq:autosys}).
Indeed, in Corollary \ref{Yges}, the decay rate is given by
\(1-\frac{a_W^2-|\varepsilon|\bar a_V}{\bar a_V a_W}\) and take values
in \((0,1)\) as
\(|\varepsilon| \leq \varepsilon^* < \frac{a_W^2}{\bar a_V}\) as
required by Corollary \ref{Yges}, hence
\(\left(1-\frac{a_W^2-|\varepsilon|}{\bar a_V a_W}\right)^k\to0\) as
\(k\to\infty\). Furthermore, the estimated decay rate can be tuned via
\(\varepsilon\) from \(1\) to \(1-\frac{a_W}{\bar a_V}\) as
\(|\varepsilon|\) decreases to zero. We can therefore make the decay
smaller by adjusting \(c_\text{stop}\), as in Theorem \ref{Vestimates}.
Hence, by tuning \(\varepsilon\), we can tune how fast the closed-loop
converges to the attractor \(\{x:\sigma(x)=0\}\), and the price to pay
is more computations in general.

\hypertarget{policy-performance-guarantees}{%
\subsection{\texorpdfstring{Policy performance guarantees
\label{near-opti}}{Policy performance guarantees }}\label{policy-performance-guarantees}}

In Section \ref{openloopoptimality}, we have provided relationships
between the finite-horizon cost \(\V[\varepsilon]\) and the
infinite-horizon cost \(\V[\infty]\). This is an important feature of
VI, but this does not directly provide us with information on the actual
value of the cost function (\ref{eq:Jinfty}) along solutions to
(\ref{eq:autosys}). Therefore, we analyse the running cost
\citep{gruneperformance} defined as \begin{equation}
\begin{split}
        \mathcal{V}_{\varepsilon}^{\text{run}}(x) := \Bigg\{\sum_{k=0}^\infty
        \ell_{\mathcal{U}^*_{\varepsilon}(\phi(k,x))}(\phi(k,x)) : \phi(\cdot,x)\text{ is a solution to (\ref{eq:autosys})}\Bigg\},
 \label{eq:Vrun}
 \end{split}
\end{equation} where
\(\ell_{\mathcal{U}^*_{\varepsilon}(\phi(k,x))}(\phi(k,x))\) is the
actual stage cost incurred at time step \(k\). It has to be noted that
\(\mathcal{V}_{\varepsilon}^\text{run}(x)\) is a set, since solutions of
(\ref{eq:autosys}) are not necessarily unique. Each element
\(V_{\varepsilon}^{\text{run}}(x) \in \mathcal{V}_{\varepsilon}^{\text{run}}(x)\)
corresponds then to the cost of a solution of (\ref{eq:autosys}).
Clearly, \(V_\varepsilon^{\text{run}}(x)\) is not necessarily bounded,
as the stage costs may not decrease to 0 in view of Theorem
\ref{algostab}. Indeed, only practical convergence is ensured in Theorem
\ref{algostab} in general. On the other hand, when the set
\(\{x\in\mathbb{R}^n : \sigma(x)=0\}\) is globally exponentially stable
as in Corollary \ref{Yges}, the elements of
\(\mathcal{V}_{\varepsilon}^{\text{run}}(x)\) in (\ref{eq:Vrun}) are
bounded and satisfy the next property. \vspace{-0.5em}

\begin{theorem}\label{Vrunestimates}
Consider  system  (\ref{eq:autosys})   and  suppose  the conditions of Corollary  \ref{Yges}  hold.  For   any  $\varepsilon$  such  that
$|\varepsilon|<\varepsilon^*$,      $x      \in       \mathbb{R}^n$,      and      $V_{\varepsilon}^{\text{run}}(x)      \in
\mathcal{V}_{\varepsilon}^{\text{run}}(x)$, it follows that\ifnotproof{ $ \V[\infty](x)\leq V_\varepsilon^{\text{run}}(x) \leq     \V[\infty](x) + w_{\varepsilon}\sigma(x)$},
\ifproof{\begin{equation}
  \V[\infty](x)\leq V_\varepsilon^{\text{run}}(x) \leq     \V[\infty](x) + w_{\varepsilon}\sigma(x),
\label{eq:Vrunestimates}
\end{equation}}
with   $  w_{\varepsilon}:=\displaystyle\frac{\bar   a_V^3}{a_W}\frac{|\varepsilon|}{a_W^2-\bar  a_V|\varepsilon|}$, where the constants come from Corollary \ref{Yges}. $\hfill\square$
\end{theorem}

The inequality \(\V[\infty](x) \leq V^{\text{run}}_{\varepsilon}(x)\) of
Theorem \ref{Vrunestimates} directly follows from the optimality of
\(\V[\infty]\). The inequality
\(\V[\varepsilon]^{\text{run}}(x) \leq \V[\infty](x) + w_{\varepsilon}\sigma(x)\)
provides a relationship between the running cost
\(V_{\varepsilon}^\text{run}(x)\) and the infinite-horizon cost at state
\(x\), \(\V[\infty](x)\). The inequality
\(V_\varepsilon^{\text{run}}(x) \leq \V[\infty](x) + w_{\varepsilon}\sigma(x)\)
confirms the intuition coming from Theorem \ref{Vestimates} that a
smaller stopping criterion leads to tighter near-optimality guarantees.
That is, when \(|\varepsilon|\to0\), \(w_{\varepsilon}\to0\) and
\(V^{\text{run}}_{\varepsilon}(x)\to\V[\infty](x)\) for any
\(x\in\mathbb{R}^n\), provided that Corollary \ref{Yges} holds. In
contrast with Theorem \ref{Vestimates}, stability of system
(\ref{eq:autosys}) is essential in Theorem \ref{Vrunestimates}. Indeed,
the term \(\scriptstyle\frac{1}{a_W^2-\bar a_V |\varepsilon|}\) in the
expression of \(w_\varepsilon\) shows that the running cost is large
when \(|\varepsilon|\) is close to \(\frac{a_W^2}{\bar a_V}\), hence,
when stability is not guaranteed, the running cost might be unbounded.

\vspace{-0.5em}

\hypertarget{illustrative-example}{%
\section{\texorpdfstring{Illustrative Example
\label{S:example}}{Illustrative Example }}\label{illustrative-example}}

We consider the discrete cubic integrator, also seen in
\citep{grimm2005,granzotto2020}, which is given by
\ifnotproof{$(x_1^+,x^+_2)= (x_1+u,x_2+u^3)$}\ifproof{\begin{equation}
\begin{split}
x_1^+&=x_1+u\\
 x_2^+&=x_2+u^3,  \label{eq:cubicsys}
 \end{split}
\end{equation}} where
\((x_1,x_2):=x\in\mathbb{R}^2\) and \(u\in\mathbb{R}\). Let
\(\sigma(x)=|x_1|^3+|x_2|\) and consider cost (\ref{eq:J}) with
\(\ell(x,u)=|x_1|^3+|x_2|+|u|^3\) for any
\((x,u)\in\mathbb{R}^2\times\mathbb{R}\). It is shown in
\citep{granzotto2020} that SA\ref{SAa} holds with
\(\overline{\alpha}_V=14\mathbb{I}\) and \(\alpha_W:=\mathbb{I}\).

Because it is notoriously difficult to exactly compute \(\V[d](x)\) and
associated sequence of optimal inputs for every \(x\in\mathbb{R}^2\), we
use an approximate scheme. In particular, we rely on a simple finite
difference approximation, with \(N=340^2\) points equally distributed in
\([-10,10]\times[-10^3,10^3]\) for the state space or, equivalently,
\(\{x\in\mathbb{R}^n : \sigma(x) \leq 2000\}\), and \(909\) equally
distributed quantized inputs in \([-20,20]\) centered at 0. We consider
three types of stopping criteria for which \(\varepsilon\) is a scalar.
For each stopping criterion, we discuss the type of guaranteed stability
and we provide in Table \ref{ExampleSTOPtable} the corresponding horizon
for different values of \(\varepsilon\), which is related to the
computation cost. Then, for each horizon, we give in Table
\ref{ExampleVIqualitable} estimates of the running cost for initial
condition \(x=(10,-10^3)\), by computing the sum in \eqref{eq:Vrun} up
to \(k=40\) instead of \(k=\infty\), as well as the value of
\(\sigma(\phi(40,x))\) to evaluate the convergence accuracy of the
corresponding policy.

We first take the uniform stopping criterion uniform stopping criterion
as in (\ref{eq:uniformliterature}), like in,
e.g.,\citep{Bertsekas:12,Sutton,pang2019adaptive,kiumarsi2017h,derongliu2015},
i.e.~\(c_\text{stop}(\varepsilon,x):=|\varepsilon|,\) with different
values of \(\varepsilon\). In this case, we have no global exponential
stability or performance guarantees like in Corollary \ref{Yges} and
Theorem \ref{Vrunestimates} a priori. Only near-optimal guarantees as in
Theorem \ref{Vestimates} and semiglobal practical stability as in
Theorem \ref{algostab} hold. For instance, by taking
\(\varepsilon=0.01\), Theorem \ref{Vestimates} holds with
\(v_\varepsilon(x) = 14 \cdot 0.01 = 0.14\) for any
\(x\in\mathbb{R}^n\).

We also consider the following relative stopping criterion, for any
\(x\in\mathbb{R}^n\) and \(\varepsilon\in\mathbb{R}\),
\(c_\text{stop}(\varepsilon,x):=|\varepsilon| \sigma(x)\). The
exponential stability of Corollary \ref{Yges} holds for any
\(\varepsilon\in\mathbb{R}\) such that
\(|\varepsilon|<\frac{a_W^2}{\bar a_V} = \frac{1}{14}\) in this case.
Moreover, we have near-optimality and performance properties as in
Theorems \ref{Vestimates} and \ref{Vrunestimates}, which were not
available for the previous stopping criterion
\(c_\text{stop}(\varepsilon,x)=|\varepsilon|\). Moreover, for
\(\varepsilon=0.01<\frac{1}{14}\), Theorem \ref{Vestimates} holds with
\(v_\varepsilon(x) = 14 \cdot 0.01 \cdot \sigma(x) = 0.14 \sigma(x)\)
for any \(x\in\mathbb{R}^n\), which is small when \(\sigma(x)\) is
small, and vice versa. Compared to the previous stopping criterion,
which leads to constant guaranteed near-optimality bound, here we have
better guarantees when \(\sigma(x)\) is small (and worse ones when
\(\sigma(x)\) is large). We observe less computations for better a
priori near-optimality properties for states near the attractor,
i.e.~when \(\sigma(x)<1\), when compared to the previous stopping
criterion. We finally consider the mixed stopping criterion
\(c_\text{stop}(\varepsilon,x):=|\varepsilon| \min\{\sigma(x),1\}\),
which provides better near-optimality guarantees than both considered
stopping criteria. We see from Table \ref{ExampleSTOPtable}, and Table
\ref{ExampleVIqualitable}, that by increasing iterations, we usually
obtain smaller and thus better running costs as well as tighter
convergence properties.

Compared to previous work \citep{granzotto2020}, where stability
properties to (approximate) value iteration are given, we require a
smaller number of iterations. Indeed, in view of
\citep[Corollary 2]{granzotto2020},
\(d\geq\bar d= \floor*{\frac{0-\ln 14^2}{\ln 13 -\ln 14}}= 71\). Of
course, this analysis is conservative and a different derivation of
\(\overline{\alpha}_V\) might provide different bounds on \(\bar d\).
Here, as the algorithm is free to choose the required number of
iterations via the stopping criterion, we significantly reduce its
conservatism. This induces smaller computational complexity, as,
e.g.~for \(\varepsilon=0.01<\frac{1}{14}\), exponential stability is
ensured with the stoping criterion verified at \(d=6\), that is,
\(8.5\%\) of iterations required by the lower bound \(\bar d =71\) of
\citep{granzotto2020}.

\begin{table}
\small
\begin{center}
\begin{tabular}{
 r l
 | c c   c c c c c
}
\bottomrule
                                                                           & & \multicolumn{7}{c}{$\bm{\varepsilon}$}  \\[-0.1em]
                                                                                       &   &  10 & 0.75 & 0.1 & 0.075& 0.05& 0.025 & 0.005    \\ 
\bottomrule                                                        
  \multirow{3}{*}{$\bm{c_\text{stop}(\varepsilon,x):}$}     & $\bm{|\varepsilon|}$& $d=6$& $d=7$& $d=8$& $d=8$& $d=8$ &$d=8$  & $d=9$       \\
                                             & $\bm{|\varepsilon|\sigma(x)}$ & $d=0$ & $d=1$& $d=3$& $d=4$& $d=5$ &$d=6$ & $d=7$       \\
                                & $\bm{|\varepsilon|\min\{\sigma(x),1\}}$ & $d=6$& $d=7$& $d=8$& $d=8$& $d=8$ &$d=8$  & $d=9$       \\\bottomrule
\end{tabular}
\caption{Required iterations to fulfill  each stopping criteria  for $N=340^2$ points equally distributed  in $\{z\in\Reals^n : \sigma(z) \leq 2000\}$. \label{ExampleSTOPtable}}
\vspace{-2em}
\end{center}
\normalsize
\end{table}
\begin{table}
\small
\begin{center}
\begin{tabular}{
 c 
 | c c   c c c c c c c
}
\bottomrule     
                                                    &  \multicolumn{9}{c}{$\bm{d}$}  \\[-0.1em]     
                                   &  0  &  1 & 3 & 4 & 5 & 6 & 7 & 8 & 9 \\ 
\bottomrule                                                        
        $\bm{V^\text{run}_d(x)}$       &  77313  &  45497 &  19931 & 19965 &  19802 & 20090 & 20359 & 20261  & 20261 \\ 
        $\bm{\sigma(\phi(40,x))}$      &  1982  &  1138 &  2.56  & 2.84 & 1.71 & 2.25 & 1.84 & 1.62 & 1.62 \\        \bottomrule
\end{tabular}
\caption{Estimation of the running cost  $V^\text{run}_d(x)$ for  $x=(10,-10^3)$ and the value of $\sigma(\phi(40,x))$. \label{ExampleVIqualitable}}
\vspace{-2em}
\end{center}
\normalsize
\end{table}

\vspace{-1em}

\hypertarget{concluding-remarks}{%
\section{\texorpdfstring{Concluding remarks
\label{conclusion}}{Concluding remarks }}\label{concluding-remarks}}

Future work includes relaxing the initial condition for VI and the main
assumptions. Another direction is extending the work towards stochastic
problems and online algorithms, towards the final goal of
stability-based computational-performance tradeoffs in reinforcement
learning.

\ifproofs
\section{Proofs}\label{sec:proofs}

\subsection{Proof of Proposition \ref{prop:terminal} }

Let \(x\in\mathbb{R}^n\) and \(d\in\mathbb{Z}_{>0}\). Since
\(\V[d-1](x)\) is the optimal value function associated to the
\((d-1)\)-horizon cost (\ref{eq:Vd}), it follows
\(\V[d-1](x)\leq \sum\limits_{k=0}^{d-1}\ell^*_d(k,x)\), and by
definition of \(\ell_d^*\),
\(\V[d](x)= \sum\limits_{k=0}^{d}\ell^*_d(k,x)\). Hence
\(\ell^*_d(d,x)=\sum\limits_{k=0}^{d}\ell^*_d(k,x)-\sum\limits_{k=0}^{d-1}\ell^*_d(k,x) \leq \V[d](x)- \V[d-1](x)\),
which gives the desired result.

\subsection{Proof of Theorem \ref{theo:terminates}}

We first consider the case when item (i) of Assumption \ref{a:stop}
holds. Let \(\Delta>0\) and \(x\in\mathbb{R}^n\) such that
\(\sigma(x)\leq \Delta\), \(\varepsilon\in \mathbb{R}^{n_\varepsilon}\)
and \(\underline\epsilon>0\) as in item (i) of Assumption \ref{a:stop}.
In view of \citep[Theorem 3]{granzotto2020} with \(\gamma=1\), it
follows that
\begin{equation}\V[d](x) \leq \V[\infty](x) \leq \V[d](x) + \overline{\alpha}_V\circ \underline{\alpha}_Y^{-1}\circ \left(\mathbb{I}-{\alpha}_Y\circ\overline{\alpha}_Y^{-1}\right)^{d}\circ\overline{\alpha}_Y(\sigma(x)),\end{equation}
for all \(d\in\mathbb{Z}_{>0}\) and \(x\in\mathbb{R}^n\), where
\(\overline{\alpha}_Y:=\overline{\alpha}_V\) and
\(\underline{\alpha}_Y,{\alpha}_Y:=\alpha_W\). Hence
\(-\V[d-1](x) - \overline{\alpha}_V\circ \underline{\alpha}_Y^{-1}\circ \left(\mathbb{I}-{\alpha}_Y\circ\overline{\alpha}_Y^{-1}\right)^{d-1}\circ\overline{\alpha}_Y(\sigma(x))\leq -\V[\infty](x)\)
and since \(\V[d](x)\leq \V[\infty](x)\), we derive that
\begin{equation}\V[d](x) -\V[d-1](x)  - \overline{\alpha}_V\circ \underline{\alpha}_Y^{-1}\circ \left(\mathbb{I}-{\alpha}_Y\circ\overline{\alpha}_Y^{-1}\right)^{d-1}\circ\overline{\alpha}_Y(\sigma(x)) \leq \V[\infty](x)-\V[\infty](x)=0,\end{equation}
thus
\begin{equation}\label{eq:diff} \V[d](x) -\V[d-1](x)   \leq \overline{\alpha}_V\circ \underline{\alpha}_Y^{-1}\circ \left(\mathbb{I}-{\alpha}_Y\circ\overline{\alpha}_Y^{-1}\right)^{d-1}\circ\overline{\alpha}_Y(\sigma(x)).\end{equation}
When \(\sigma(x)=0\),
\(\V[d](x)-\V[d-1](x)\leq\overline{\alpha}_V(0)=0\). Hence,
\(\V[d](x)- \V[d-1](x)\leq c_\text{stop}(\varepsilon,x)\) for any
\(d\in\mathbb{Z}_{>0}\) since \(c_{\text{stop}}(\varepsilon,x)\geq0\).
When \(\sigma(x)>0\),
\(\overline{\alpha}_V\circ \underline{\alpha}_Y^{-1}\circ \left(\mathbb{I}-{\alpha}_Y\circ\overline{\alpha}_Y^{-1}\right)^{d-1}\circ\overline{\alpha}_Y(\sigma(x))\)
is upper-bounded by
\(\overline{\alpha}_V\circ \underline{\alpha}_Y^{-1}\circ \left(\mathbb{I}-{\alpha}_Y\circ\overline{\alpha}_Y^{-1}\right)^{d-1}\circ\overline{\alpha}_Y(\sigma(\Delta))\)
as \(\sigma(x)\leq \Delta\) since the considered function is
non-decreasing. The term
\(\overline{\alpha}_V\circ \underline{\alpha}_Y^{-1}\circ \left(\mathbb{I}-{\alpha}_Y\circ\overline{\alpha}_Y^{-1}\right)^{d-1}\circ\overline{\alpha}_Y(\sigma(\Delta))\)
can be made arbitrarily close to \(0\) by increasing \(d\), according to
\citep[Lemma  3]{granzotto2020} as \(\gamma=1\) and
\(\sigma(x)\leq\Delta\). In particular, we take \(d^*\) sufficiently
large such that
\(\overline{\alpha}_V\circ \underline{\alpha}_Y^{-1}\circ \left(\mathbb{I}-{\alpha}_Y\circ\overline{\alpha}_Y^{-1}\right)^{d^*-1}\circ\overline{\alpha}_Y(\sigma(\Delta))\leq\underline\epsilon\)
where \(\underline\epsilon\) comes from item (i) of Assumption
\ref{a:stop} and is such that
\(\underline\epsilon\leq c_{\text{stop}}(\varepsilon,x)\). Hence
\(\V[d](x) -\V[d-1](x)\leq\overline{\alpha}_V\circ \underline{\alpha}_Y^{-1}\circ \left(\mathbb{I}-{\alpha}_Y\circ\overline{\alpha}_Y^{-1}\right)^{d-1}\circ\overline{\alpha}_Y(\sigma(x))\leq c_{\text{stop}}(\varepsilon,x)\)
for any \(d\geq d^*\) and \(d^*\) depends on \(\Delta\) a priori, but
not on \(x\). We have proved the desired result.

Consider now the case when item (ii) of Assumption \ref{a:stop} holds
with \(L<\infty\). Let \(\Delta>0\) and \(x\in\mathbb{R}^n\) such that
\(\sigma(x)\leq \Delta\), \(\varepsilon\in \mathbb{R}^{n_\varepsilon}\)
and \(\underline\epsilon>0\) as in item (ii) of Assumption \ref{a:stop}
and is such that
\(\underline\epsilon \sigma(x)\leq c_\text{stop}(\varepsilon,x)\). It
follows from the inequalities of item (ii) of Assumption \ref{a:stop}
that there exists some \(\delta>0\) such that
\(\overline{\alpha}_V\circ \underline{\alpha}_Y^{-1}\circ \left(\mathbb{I}-{\alpha}_Y\circ\overline{\alpha}_Y^{-1}\right)^{d}\circ\overline{\alpha}_Y(s)\leq \frac{\bar a_V^2}{a_W} \left(1-\frac{a_W}{\bar a_V}\right)^{d} s\)
for any \(s\in[0,\delta]\), see
\citep[proof of Corollary 1, equation (47)]{granzotto2020}. Note that
\(1-\frac{a_W}{\bar a_V}\in(0,1)\), hence
\(\left(1-\frac{a_W}{\bar a_V}\right)^{d-1}\) can be made as small as
desired. Therefore, there exists \(\bar d\) such that
\(\frac{\bar a_V^2}{a_W} \left(1-\frac{a_W}{\bar a_V}\right)^{\bar d-1} \leq \underline\epsilon\)
where \(\underline\epsilon\). Thus, when \(\sigma(x)\in[0,\delta]\) and
in view of (\ref{eq:diff}), it follows that
\(\V[d](x) -\V[d-1](x) \leq \overline{\alpha}_V\circ \underline{\alpha}_Y^{-1}\circ \left(\mathbb{I}-{\alpha}_Y\circ\overline{\alpha}_Y^{-1}\right)^{d}\circ\overline{\alpha}_Y(\sigma(x)) \leq \frac{\bar a_V^2}{a_W} \left(1-\frac{a_W}{\bar a_V}\right)^{d} \sigma(x) \leq \underline\epsilon \sigma(x) \leq c_\text{stop}(\varepsilon,x)\)
for any \(d\geq \bar d\). When \(\sigma(x)\in[\delta,\Delta]\), we have
\(c_\text{stop}(\varepsilon,x) \geq \underline\epsilon \delta>0\), and
we recover \(d^*\) such that
\(\V[d](x) -\V[d-1](x)\leq c_{\text{stop}}(\varepsilon,x)\) for
\(\sigma(x)\in[\delta,\Delta]\) by applying the steps made above, for
when item (i) of Assumption \ref{a:stop} holds. Therefore, it follows
that, for any \(d\geq\min\{\bar d, d^*\}\),
\(\V[d](x) -\V[d-1](x)\leq c_{\text{stop}}(\varepsilon,x)\) holds for
any \(\sigma(x)\in[0,\Delta]\). We have proved the desired result.

Suppose now that item (ii) of Assumption \ref{a:stop} holds with
\(L=\infty\). Let \(x\in\mathbb{R}^n\),
\(\varepsilon\in \mathbb{R}^{n_\varepsilon}\) and
\(\underline\epsilon>0\) as in item (ii) of Assumption \ref{a:stop}. It
follows then that
\(\overline{\alpha}_V\circ \underline{\alpha}_Y^{-1}\circ \left(\mathbb{I}-{\alpha}_Y\circ\overline{\alpha}_Y^{-1}\right)^{d}\circ\overline{\alpha}_Y(\sigma(x))\leq \frac{\bar a_V^2}{a_W} \left(1-\frac{a_W}{\bar a_V}\right)^{d} \sigma(x)\)
for all \(x\in\mathbb{R}^n\). Hence, by invoking the same arguments as
above when item (ii) of Assumption \ref{a:stop} holds, there exists
\(\bar d>0\) such that
\(\V[d](x)- \V[d-1](x)\leq c_\text{stop}(\varepsilon,x)\) for any
\(d\geq\bar d\). The proof is complete.

\subsection{Proof of Theorem \ref{Vestimates}}

Let \(x\in\mathbb{R}^n\) such that \(\sigma(x)\leq\Delta\) and
\(\varepsilon\in\mathbb{R}^{n_\varepsilon}\),
\(d\in\mathbb{Z}_{\geq 0}\) as in (\ref{eq:stop}), which exists since
Theorem \ref{theo:terminates} holds. Hence, optimal sequence
\([\us[0],\us[1],\ldots,\us[d]]:=\ustar[\varepsilon](x)\) and cost
\(\V[\varepsilon](x)\) defined in (\ref{eq:V}) are well-defined. Since
\(\V[\varepsilon]\) is a finite-horizon optimal cost,
\(\V[\varepsilon](x)\leq\V[\infty](x)\). On the other hand, consider the
infinite-horizon sequence
\(\bm{u}=[\us[0],\us[1],\ldots \us[d-1],\ustar[\infty](\phi(d,x,\ustar[\varepsilon](x)|_{d}))]\),
where \(\ustar[\infty](\phi(d,x,\ustar[\varepsilon](x)|_{d}))\) exists
in view of item (i) of SA\ref{SAa}. It follows from the optimality of
\(\V[\infty](x)\) that \(\V[\infty](x)\leq\J[\infty](x,\bm{u})\), and
from the definition of \(\bm{u}\) that
\(\J[\infty](x,\bm{u})= \V[\varepsilon](x)+\V[\infty](\phi(d,x,\ustar[\varepsilon](x)|_{d}))\),
which is finite. By invoking item (i) of SA\ref{SAa}, we derive
\(\V[\infty](x)\leq\V[\varepsilon](x)+\overline{\alpha}_V(\sigma(\phi(d,x,\ustar[\varepsilon](x)|_{d})))\).
In view of Lemma \ref{prop:terminal} and item (ii) of SA\ref{SAa},
\(\sigma(\phi(d,x,\ustar[\varepsilon](x)|_{d}))\leq\alpha_W^{-1}(c_{\text{stop}}(\varepsilon,x))\),
thus
\(\V[\infty](x)\leq\V[\varepsilon](x)+ \overline{\alpha}_V\circ\alpha_W^{-1}(c_{\text{stop}}(\varepsilon,x))\)
and the proof is complete.

\subsection{Proof of Theorem \ref{algostab}}

First, we prove the following result, which provides Lyapunov properties
for system (\ref{eq:autosys}) that we use to derive the main stability
result afterwards.

\begin{theorem}\label{YLyapunovProp}
Let $\Y[]:=\V[\infty]$, the following holds.
\begin{itemize}[itemindent=\widthof{(ii)}]
   \item[(i)]  For any $x\in\mathbb{R}^n$, \begin{equation}\underline{\alpha}_Y(\sigma(x)) \leq \Y[](x)  \leq \overline{\alpha}_Y(\sigma(x)),\end{equation}
where $\underline{\alpha}_Y:=\alpha_W,\overline{\alpha}_Y:=\overline{\alpha}_V$, with $\alpha_W,\overline{\alpha}_V$ from SA\ref{SAa}.
   \item[(ii)] For  any $x\in\mathbb{R}^n$, $\varepsilon\in\mathbb{R}^{n_\varepsilon}$,  $v \in F^*_{\varepsilon}(x)$,  \begin{equation} \Y[](v)-\Y[](x)
   \leq -{\alpha}_Y(\sigma(x)) + \overline{\alpha}_V\circ\alpha_W^{-1}(c_\text{stop}(\varepsilon,x))\end{equation}
     where ${\alpha}_Y=\alpha_W$, with $\alpha_W$ and $\overline{\alpha}_V$ from SA\ref{SAa}, and $c_\text{stop}$ comes from (\ref{eq:stop}). $\hfill\square$
\end{itemize}
\end{theorem}

\noindent\textbf{Proof.} Let
\(\varepsilon\in\mathbb{R}^{n_\varepsilon}\), \(x\in\mathbb{R}^n\) and
\(v \in F^* _ {\varepsilon}(x)\), which is well-defined in view of
Theorem \ref{theo:terminates}. There exists
\([\us[0],\us[1],\ldots,\us[d]]=\bm{u}^*_{\varepsilon}(x)\) such that
\(v=f(x,\us[0])\) and \(\bm{u}^*_{\varepsilon}(x)\) is an optimal input
sequence for system (\ref{eq:sys}) and cost (\ref{eq:J}) with horizon
\(d\), hence \(\V[\varepsilon](x)=\J[d](x,\bm{u}^*_{\varepsilon}(x))\).
Moreover, in view of Lemma \ref{prop:terminal} and item (ii) of
SA\ref{SAa},
\(\sigma(\phi(d,x,\ustar[\varepsilon](x)|_{d}))\leq\alpha_W^{-1}(c_{\text{stop}}(\varepsilon,x))\).

From items (i) and (ii) of SA\ref{SAa}, we have
\(\Y[](x)=\V[\infty](x)\leq\overline{\alpha}_V(\sigma(x))=:\overline{\alpha}_Y(\sigma(x))\).
On the other hand, we have from item (ii) of SA\ref{SAa} that
\(\alpha_W(\sigma(x))\leq \ell^*_d(0,x)\). This implies that
\(\alpha_W(\sigma(x))\leq \V[\varepsilon](x)\leq\V[\infty](x)=\Y[](x)\).
Hence item (i) of Theorem \ref{YLyapunovProp} holds with
\(\underline{\alpha}_Y=\alpha_W\).

Consider the sequence
\(\hat{\bm{u}}:=[\us[1],\us[2],\ldots,\us[d-1],\bar{\bm{u}}]\) where
\(\bar{\bm{u}}:=\bm{u}^*_{\infty}( \phi(d , x , \bm{u}^*_{\varepsilon}(x)|_{d}) )\),
\(\bm{u}^*_{\varepsilon}(x)|_{d}=[\us[0],\ldots,\us[d-1]]\) and \(\phi\)
denotes the solution of system (\ref{eq:sys}). The sequence
\(\hat{\bm{u}}\) consists of the first \(d\) elements of
\(\bm{u}^*_{\varepsilon}(x)\) after \(\us[0]\), followed by an optimal
input sequence of infinite length at state
\(\phi(d,x,\ustar[\varepsilon](x)|_{d})\), which exists according to
item (i) of SA\ref{SAa}. Sequence \(\bar{\bm{u}}\) minimizes
\(\J[\infty](\phi(d,x,\ustar[\varepsilon](x)|_{d}),\bar{\bm{u}})\) by
virtue of item (i) of SA\ref{SAa}. From the definition of cost \(\J\) in
(\ref{eq:J}) and \(\V[\infty](v)\) in view of item (i) of SA\ref{SAa},
\begin{equation}
\label{eq:propeqVa}
\begin{split}
\V[\infty](v)\quad
  &\leq\quad
  \J[\infty](v,\hat{\bm{u}})\\
  &=\quad
     \J[d-1](v,\hat{\bm{u}}|_{d-1})\\&\quad+\J[\infty](\phi(d-1,v,\hat{\bm{u}}|_{d-1}),\bar{\bm{u}}).
     \end{split}
\end{equation} From Bellman optimality principle, we have
\(\V[\varepsilon](x) = \V[d](x)=\ell^*_d(0,x) + \V[d-1](v) = \ell^*_d(0,x)+\J[d-1](v,\hat{\bm{u}}|_{d-1})\),
hence
\begin{equation}\J[d-1](v,\hat{\bm{u}}|_{d-1})=V_{\varepsilon}(x)-\ell^*_d(0,x)\label{eq:propeqVb}.
\end{equation} Moreover, by item (i) of SA\ref{SAa}, \begin{align*}
\MoveEqLeft \J[\infty](\phi(d-1,v,\hat{\bm{u}}|_{d-1}),\bar{\bm{u}}) \\&\leq\overline{\alpha}_V(\sigma(\phi(d-1,v,\hat{\bm{u}}|_{d-1}))).\addtocounter{equation}{1}\tag{\theequation}\label{eq:propeqVc}
\end{align*} Consequently, in view of (\ref{eq:propeqVa}),
(\ref{eq:propeqVb}) and (\ref{eq:propeqVc}), \begin{align*}
\MoveEqLeft \V[\infty](v)  \leq \V[\varepsilon](x) -\ell^*_d(0,x)\\
     &\qquad+\overline{\alpha}_V(\sigma(\phi(d-1,v,\hat{\bm{u}}|_{d-1}))).
\addtocounter{equation}{1}\tag{\theequation}\label{eq:propeqa}
\end{align*} Since
\(\phi(d-1,v,\hat{\bm{u}}|_{d-1})=\phi(d,x,\ustar[\varepsilon](x)|_{d})\)
and
\(\sigma(\phi(d,x,\ustar[\varepsilon](x)|_{d}))\leq\alpha_W^{-1}(c_{\text{stop}}(\varepsilon,x))\)
holds, it follows \begin{equation}
\V[\infty](v)
     \leq \V[\varepsilon](x)-\ell^*_d(0,x)+\overline{\alpha}_V\circ\alpha_W^{-1}(c_\text{stop}(\varepsilon,x)).
\label{eq:propeqb}
\end{equation} By Theorem \ref{Vestimates},
\(\V[\varepsilon](x)\leq\V[\infty](x)\), thus \begin{equation}
\V[\infty](v)
     \leq \V[\infty](x)-\ell^*_d(0,x)+\overline{\alpha}_V\circ\alpha_W^{-1}(c_\text{stop}(\varepsilon,x)).
\label{eq:propeqfinal}
\end{equation} By invoking item (ii) of SA\ref{SAa}, we derive
\(\V[\infty](v) \leq \V[\infty](x) -\alpha_W(\sigma(x))+\overline{\alpha}_V\circ\alpha_W^{-1}(c_\text{stop}(\varepsilon,x))\),
and since \(\Y[]=\V[\infty]\), the proof is completed with
\({\alpha}_Y:=\alpha_W\). \(\hfill\blacksquare\)

\vspace{1em}

Item (i) states that \(\Y[]\) is positive definite and radially
unbounded with respect to the set \(\{x:\sigma(x)=0\}\). Item (ii) of
Theorem \ref{YLyapunovProp} shows that \(\Y[]\) strictly decreases along
the solutions to (\ref{eq:autosys}) up to a perturbative term
\(\overline{\alpha}_V\circ\alpha_W^{-1}(c_\text{stop}(\varepsilon,x))\),
which can be made as small as desired by selecting \(|\varepsilon|\)
close to \(0\) as
\(\overline{\alpha}_V\circ\alpha_W^{-1}(c_\text{stop}(\varepsilon,x))\leq\overline{\alpha}_V\circ\alpha_W^{-1}(\theta(|\varepsilon|,\sigma(x)))\),
per Assumption \ref{cstop2}.

We are ready to prove Theorem \ref{algostab}. The proof heavily borrows
from \citep[Theorem2]{granzottoCDC2019} and
\citep[Theorem 3.3]{granz2019}, we nevertheless include it in this
technical report for self-completeness. Let \(\Delta,\delta>0\). We
select \(\varepsilon^*>0\) such that \begin{alignat*}{3}
   \theta(\varepsilon^*,\underline{\alpha}_Y^{-1}(\widetilde\Delta))  &< \alpha_W\circ\overline{\alpha}_V^{-1}(\frac{1}{2}\widetilde{\alpha}_Y(\widetilde\delta)), \addtocounter{equation}{1}\tag{\theequation}\label{eq:gdcondC}
\end{alignat*} where
\(\widetilde{\alpha}_Y:=\alpha_W\circ\overline{\alpha}_Y^{-1}\),
\(\widetilde\Delta:= \overline{\alpha}_Y(\Delta)\),
\(\widetilde\delta:= \left(\mathbb{I}-\frac{\widetilde{\alpha}_Y}{2}\right)^{-1}\circ\underline{\alpha}_Y(\delta)\)
and \(\theta\) comes from Assumption \ref{cstop2}. Note that
\(\left(\mathbb{I}-\frac{\widetilde{\alpha}_Y}{2}\right)^{-1}\) is
indeed of class \(\mathcal{K}_\infty\) as we assume without loss of
generality that\footnote{See \citep{granz2019} for details.}
\(\mathbb{I}-\widetilde{\alpha}_Y\in\mathcal{K}_\infty\), hence
\(\mathbb{I}-\widetilde{\alpha}_Y+\frac{\widetilde{\alpha}_Y}{2}\in\mathcal{K}_\infty\)
and so is its inverse. Inequality (\ref{eq:gdcondD}) can always be
verified by taking \(\varepsilon^*\) sufficiently small since
\(\theta(\cdot,\overline{\alpha}_Y^{-1}(\widetilde\Delta))\in\mathcal{K}\),
and
\(\alpha_W\circ\overline{\alpha}_V^{-1}(\frac{1}{2}\widetilde{\alpha}_Y(\widetilde\delta))>0\).
It follows from \(\theta(\cdot,s)\in\mathcal{K}\) for any \(s>0\) and
\(\theta(s,\cdot)\) is non-decreasing for any \(s\geq0\), that
\(\theta(|\varepsilon|,\underline{\alpha}_Y^{-1}(s))\leq\theta(\varepsilon^*,\underline{\alpha}_Y^{-1}(\widetilde\Delta))\)
for any \(s\in[0,\widetilde\Delta]\) and
\(|\varepsilon|<\varepsilon^*\). Furthermore, from Assumption
\ref{cstop2} and item (i) of Theorem \ref{YLyapunovProp}, we derive
\(c_\text{stop}(\varepsilon,x)\leq \theta(|\varepsilon|,\underline{\alpha}_Y^{-1}(\Y[](x)))\).
Thus, in view of (\ref{eq:gdcondC}),
\begin{equation}\label{eq:cstopcond}
  c_\text{stop}(\varepsilon,x)\leq \alpha_W\circ\overline{\alpha}_V^{-1}\left(\frac{1}{2}\widetilde{\alpha}_Y(\widetilde\delta)\right)
\end{equation} for any \(x\) such that \(\Y[](x)\leq\widetilde\Delta\).
On the other hand, we have
\(\alpha_W\circ\overline{\alpha}_V^{-1}(\frac{1}{2}\widetilde{\alpha}_Y(\widetilde\delta))\leq\alpha_W\circ\overline{\alpha}_V^{-1}(\frac{1}{2}\widetilde{\alpha}_Y(s))\)
for any \(s\in[\widetilde\delta,\infty)\). Hence, for any
\(x\in\mathbb{R}^n\) such that
\(\Y[](x)\in[\widetilde\delta,\widetilde\Delta]\) and
\(|\varepsilon|<\varepsilon^*\), \begin{alignat*}{3}
   \overline{\alpha}_V\circ\alpha_W^{-1}(c_\text{stop}(\varepsilon,x))    &\leq  \frac{\widetilde{\alpha}_Y(\widetilde\delta)}{2} \leq \frac{\widetilde{\alpha}_Y(\Y[](x))}{2}.  \addtocounter{equation}{1}\tag{\theequation}\label{eq:gdcondD}
\end{alignat*}

Let \(x\in\mathbb{R}^n\) with \(\sigma(x)\leq\Delta\) and
\(v \in F^* _ {\varepsilon}(x)\). In view of (\ref{eq:cstopcond}) and
items (i) and (ii) of Theorem \ref{YLyapunovProp}, \begin{equation}
\Y[](v)-\Y[](x) \leq -\widetilde{\alpha}_Y(\Y[](x)) +\overline{\alpha}_V\circ\alpha_W^{-1}(c_\text{stop}(\varepsilon,x)).
     \label{eq:Ycontracttildepre}
\end{equation} Since \(\sigma(x) \leq \Delta\),
\(\Y[](x) \leq \overline{\alpha}_Y(\sigma(x)) \leq \overline{\alpha}_Y(\Delta) = \widetilde\Delta\).
Consider \(\Y[](x) \in [0,\widetilde\delta)\). Since
\(c_\text{stop}(\varepsilon,x)\leq \alpha_W\circ\overline{\alpha}_V^{-1}(\frac{1}{2}\widetilde{\alpha}_Y(\widetilde\delta))\)
holds for \(\Y[](x)\leq\widetilde\Delta\), it holds here. Furthermore,
since \(\mathbb{I}-\widetilde{\alpha}_Y\in\mathcal{K}_\infty\) holds
without loss of generality, and in view of (\ref{eq:Ycontracttildepre}),
\begin{equation}
\begin{split}
  \Y[](v)       &\leq   \Y[](x)-\widetilde{\alpha}_Y(\Y[](x)) + \overline{\alpha}_V(\theta(\varepsilon^*,\sigma(x)))\\
                      &\leq   \left(\mathbb{I}-\widetilde{\alpha}_Y\right)(\widetilde\delta) + \frac{1}{2}\widetilde{\alpha}_Y(\widetilde\delta).
                      \end{split}
\end{equation} Given the definition of \(\widetilde\delta\),
\begin{equation}
  \Y[](v) \leq  \left(\mathbb{I}-\frac{\widetilde{\alpha}_Y}{2}\right)(\widetilde\delta)  = \underline{\alpha}_Y(\delta). \label{eq:Yattracted}
\end{equation}

When \(\Y[](x)\geq\widetilde\delta\), we derive from (\ref{eq:gdcondD})
that
\(-\widetilde{\alpha}_Y(\Y[](x))+\overline{\alpha}_V(c_\text{stop}(\varepsilon,x))\leq-\tfrac{1}{2}\widetilde{\alpha}_Y(\Y[](x))\).
Thus, from (\ref{eq:Ycontracttildepre}), \begin{align*}
 \Y[](v)-\Y[](x) \leq -\tfrac{1}{2}\widetilde{\alpha}_Y(\Y[](x)).
\addtocounter{equation}{1}\tag{\theequation}\label{eq:Ycontracting}
\end{align*} In view of (\ref{eq:Yattracted}) and
(\ref{eq:Ycontracting}), it follows for any \(k\in\mathbb{Z}_{\geq 0}\)
that \begin{equation} \label{eq:Ysolonestep}
  \Y[](\phi(k+1,x))\leq \max\left\{(\mathbb{I}-\tfrac{1}{2}\widetilde{\alpha}_Y)(\Y[](x)), \underline{\alpha}_Y(\delta)\right\},
\end{equation} where \(\phi(k,x)\) is a solution starting at \(x\) for
system (\ref{eq:autosys}). Furthermore, when
\(\Y[](x)\leq \underline{\alpha}_Y(\delta)\),
\(\Y[](v) \leq \underline{\alpha}_Y(\delta)\) follows. Indeed, if
\(\Y[](x) \in [\widetilde\delta,\widetilde\Delta]\),
\(\Y[](v) \leq \Y[](x)\leq \underline{\alpha}_Y(\delta)\) from
(\ref{eq:Ycontracting}), and if \(\Y[](x) \in[0,\widetilde\delta)\), we
deduce \(\Y[](v) \leq \underline{\alpha}_Y(\delta)\) from
(\ref{eq:Yattracted}). Hence the set
\(\{z\in\mathbb{R}^n \,:\, \Y[](z)\leq \underline{\alpha}_Y(\delta)\}\)
is forward invariant for system (\ref{eq:autosys}). By iterating
(\ref{eq:Ysolonestep}), we obtain \begin{equation}
  \Y[](\phi(k,x)) \leq \max \left\{ \widetilde\beta(\Y[](x),k),\underline{\alpha}_Y(\delta) \right\},
\end{equation} where
\(\widetilde\beta(s,k)=\left(\mathbb{I}-\frac{1}{2}\widetilde{\alpha}_Y\right)^{(k)}(s)\)
for any \(s\geq0\), with \(\widetilde\beta\in\mathcal{KL}\) as
\(\lim_{k\to\infty}\left(\mathbb{I}-\frac{1}{2}\widetilde{\alpha}_Y\right)^{(k)}(s)=0\)
for any \(s\geq0\), since\footnote{See \citep[Lemma B.1]{granz2019}}
\(\left(\mathbb{I}-\frac{1}{2}\widetilde{\alpha}_Y\right)(s)<s\) for
\(s>0\) and
\(\left(\mathbb{I}-\frac{1}{2}\widetilde{\alpha}_Y\right)(0)=0\).
Finally, invoking
\(\underline{\alpha}_Y(\sigma(x))\leq \Y[](x) \leq \overline{\alpha}_Y(\sigma(x))\),
we deduce \begin{equation}
  \sigma(\phi(k,x)) \leq \max \left\{\underline{\alpha}_Y^{-1}\left(\widetilde\beta(\overline{\alpha}_Y(\sigma(x)),k)\right),\delta \right \}.
\end{equation} Thus Theorem \ref{algostab} holds with
\(\beta(s,k) = \underline{\alpha}_Y^{-1}\left(\widetilde\beta(\overline{\alpha}_Y(s),k)\right)\)
for any \(s\geq0\) and \(k\in\mathbb{Z}_{\geq 0}\).

\subsection{Proof of Corollary \ref{Yges}}

We follow the steps made in \citep[Corollary  3.2]{granz2019}. Let
\(x \in \mathbb{R}^n\). We select \(\varepsilon^*<\frac{a_W^2}{a_V}\) as
in Corollary \ref{Yges} and let
\(\varepsilon\in\mathbb{R}^{n_\varepsilon}\) such that
\(|\varepsilon|\leq\varepsilon^*\) and \(v \in F^*_{\varepsilon}(x)\).
In particular, from item (i) of Theorem \ref{YLyapunovProp},
\(\alpha_W(\sigma(x)\leq \Y[](x) \leq \overline{\alpha}_V(\sigma(x))\),
and since item (ii) of Assumption \ref{a:stop} holds with \(L=\infty\)
and \(a_W s \leq \alpha_W(s)\), \(\overline{\alpha}_V(s)\leq a_V s\)
with \(a_W,a_V>0\) for any \(s>0\), we obtain \begin{equation} 
a_W \sigma(x) \leq Y(x) \leq \bar a_V \sigma(x). 
\label{eq:ygesybound}
\end{equation} Similarly, in view of item (ii) of Theorem
\ref{YLyapunovProp},
\(\Y[](v)-\Y[](x) \leq -\alpha_W(\sigma(x)) +\overline{\alpha}_V\circ\alpha_W^{-1}(c_\text{stop}(\varepsilon,x)),\)
and since \(c_\text{stop}(\varepsilon,x)\leq |\varepsilon|\sigma(x)\)
holds, we derive\\
\begin{equation}
\Y[](v)-\Y[](x) \leq-\left(\frac{a_W^2-|\varepsilon|\bar
a_V}{a_W}\right)\sigma(x). 
\label{eq:ygesa} 
\end{equation} Note that for \(|\varepsilon|\leq\varepsilon^*\),
\(\left(\frac{a_W^2-|\varepsilon|\bar a_V}{a_W}\right) > 0\). Hence, in
view of \eqref{eq:ygesa} and \eqref{eq:ygesybound},
\begin{equation}\Y[](v)-\Y[](x) \leq-\left(\frac{a_W^2-|\varepsilon|\bar
a_V}{\bar  a_V a_W}\right)\Y[](x)\end{equation} holds for any
\(\sigma(x)\geq0\). Let \(x\in\mathbb{R}^n\) and denote \(\phi(k,x)\) be
a corresponding solution to (\ref{eq:autosys}) at time
\(k\in \mathbb{Z}_{\geq 0}\), it holds that
\(\Y[](\phi(k,x))\leq \left(1-\frac{a_W^2-|\varepsilon|\bar a_V}{\bar a_V a_W}\right)^{k}\Y[](x)\).
In view of \eqref{eq:ygesybound}, it follows from
\(\Y[](\phi(k,x)) \leq \left(1-\frac{a_W^2-|\varepsilon|\bar a_V}{\bar a_V a_W}\right)^{k}\Y[](x)\)
that
\(a_W \sigma(\phi(k,x))\leq \left(1-\frac{a_W^2-|\varepsilon|\bar a_V}{\bar a_V a_W}\right)^{k} \bar a_V \sigma(x)\)
hence
\(\sigma(\phi(k,x)) \leq \frac{\bar a_V}{a_W}\sigma(x)\left(1-\frac{a_W^2-|\varepsilon|\bar a_V}{\bar a_V a_W}\right)^{k}\)
and the proof is concluded.

\subsection{Proof of Theorem \ref{Vrunestimates}}

Let \(x \in \mathbb{R}^n\), \(\varepsilon\in\mathbb{R}^{n_\varepsilon}\)
such that \(|\varepsilon|<\varepsilon^*\) where \(\varepsilon^*\) is
selected as in Corollary \ref{Yges},
\(\phi(k+1,x)\in F^*_{\varepsilon}(\phi(k,x))\) for any
\(k\in\mathbb{Z}_{\geq 0}\) where \(\phi\) is a solution to
(\ref{eq:autosys}) initialized at \(x\). Let
\(u^r_k\in\mathcal{U}^*_{\varepsilon}(\phi(k,x))\) such that
\(\phi(k+1,x)=f(\phi(k,x),u^r_k)\), and note, since inputs from
(\ref{eq:autosys}) are the first input of \(\ustar[\varepsilon]\)
applied in a receding horizon fashion,
\(u^r_k=\ustar[\varepsilon](\phi(k,x))|_0\) and therefore
\(\ell(\phi(k,x),u_k^r)=\ell(\phi(k,x),\ustar[\varepsilon](\phi(k,x))|_0)= \ell^*_d(0,\phi(k,x))\)
by definition of \eqref{eq:lddef}. Consider then \begin{equation}
  V_{\varepsilon}^\text{run}(x)=  \sum_{k=0}^\infty \ell^*_d(0,\phi(k,x)),
\label{eq:Vavgdiffphi}
\end{equation} Note that indeed
\(V_{\varepsilon}^\text{run}(x)\in\mathcal{V}_{\varepsilon}^{\text{run}}(x)\).
It follows from the proof of Theorem \ref{YLyapunovProp}, in particular
\eqref{eq:propeqfinal}, that \begin{equation}
\V[\infty](\phi(k+1,x))
     \leq \V[\infty](\phi(k,x))-\ell^*_d(0,\phi(k,x))+\overline{\alpha}_V\circ\alpha_W^{-1}(c_\text{stop}(\varepsilon,\phi(k,x))),
\end{equation} from which we deduce, for any \(N\geq0\), \begin{align*}
  \MoveEqLeft \sum_{k=0}^N \ell^*_d(0,\phi(k,x))\\ 
                                &\leq\V[\infty](\phi(0,x)) -\V[\infty](\phi(1,x))+\overline{\alpha}_V\circ\alpha_W^{-1}(c_\text{stop}(\varepsilon,\phi(0,x)))\\
                               &\quad+\V[\infty](\phi(1,x)) -\V[\infty](\phi(2,x))+\overline{\alpha}_V\circ\alpha_W^{-1}(c_\text{stop}(\varepsilon,\phi(1,x)))\\
                               &\quad+\ldots\\
                               &\quad+\V[\infty](\phi(N,x)) -\V[\infty](\phi(N+1,x))\\
                               &\quad+\overline{\alpha}_V\circ\alpha_W^{-1}(c_\text{stop}(\varepsilon,\phi(N,x)))\\
                               &\leq \V[\infty](\phi(0,x))  + \sum_{k=0}^{N} \overline{\alpha}_V\circ\alpha_W^{-1}(c_\text{stop}(\varepsilon,\phi(k,x))),
\addtocounter{equation}{1}\tag{\theequation}\label{eq:Vsumdiffphi}
\end{align*} Hence, when \(N\to\infty\), \begin{equation}
\begin{split}
  V_{\varepsilon}^\text{run}(x)  &\leq \V[\varepsilon](\phi(0,x))  + \sum_{k=0}^{\infty} \overline{\alpha}_V\circ\alpha_W^{-1}(c_\text{stop}(\varepsilon,\phi(k,x))).
\end{split}\label{eq:Vsumdiffphiestimates}
\end{equation} All that remains is to compute a bound on
\(\sum_{k=0}^{\infty} \overline{\alpha}_V\circ\alpha_W^{-1}(c_\text{stop}(\varepsilon,\phi(k,x)))\),
which is possible by recalling that
\(\sigma(\phi(k,x))\leq \frac{\bar a_V}{a_W}\sigma(x)\left(1-\frac{a_W^2-|\varepsilon|\bar a_V}{\bar a_V a_W}\right)^{k}\)
holds from Corollary \ref{Yges} and
\(\overline{\alpha}_V\circ\alpha_W^{-1}(c_\text{stop}(\varepsilon,\phi(k,x)))\leq \frac{\bar a_V}{a_W} |\varepsilon|\sigma(\phi(k,x))\)
as the conditions of Corollary \ref{Yges} are assumed to hold.
Specifically,
\(\sum_{k=0}^{\infty} \overline{\alpha}_V\circ\alpha_W^{-1}(c_\text{stop}(\varepsilon,\phi(k,x))) \leq |\varepsilon|\frac{\bar a_V^2}{a_W^2}\sigma(x) \sum_{k=0}^{\infty} \left(1-\frac{a_W^2 -|\varepsilon|\bar a_V}{\bar a_V a_W}\right)^{k}\),
which provides (\ref{eq:Vrunestimates}) as
\(\sum_{k=0}^{\infty} \left(1-\frac{a_W^2-|\varepsilon|\bar a_V}{\bar a_V a_W}\right)^{k}=\frac{\bar a_V a_W}{a_W^2-\bar a_V |\varepsilon|}\).
The lower bound \(\V[\infty](x)\leq V_{\varepsilon}^\text{run}(x)\)
follows from the optimality of \(\V[\infty](x)\). Since
(\ref{eq:Vsumdiffphiestimates}) holds for an arbitrary solution of
(\ref{eq:autosys}), \(\phi(k+1,x)=f(\phi(k,x),u^r_k)\) for any
\(k\in\mathbb{Z}_{\geq 0}\), the resulting bound holds for any
\(V_{\varepsilon}^\text{run}(x)\in\mathcal{V}_{\varepsilon}^{\text{run}}(x)\).

\fi

\bibliographystyle{plain}

\bibliography{IEEEabrv,vicontrol}

\begin{thebibliography}{30}
\providecommand{\natexlab}[1]{#1}
\providecommand{\url}[1]{\texttt{#1}}
\expandafter\ifx\csname urlstyle\endcsname\relax
  \providecommand{\doi}[1]{doi: #1}\else
  \providecommand{\doi}{doi: \begingroup \urlstyle{rm}\Url}\fi

\bibitem[Anderson and Moore(2007)]{anderson2007optimal}
B.~D.~O. Anderson and J.~B. Moore.
\newblock \emph{Optimal control: linear quadratic methods}.
\newblock Courier Corporation, 2007.

\bibitem[{Arnold} and {Laub}(1984)]{ArnoldRiccati84}
W.~F. {Arnold} and A.~J. {Laub}.
\newblock Generalized eigenproblem algorithms and software for algebraic
  riccati equations.
\newblock \emph{Proceedings of the IEEE}, 72\penalty0 (12):\penalty0
  1746--1754, 1984.

\bibitem[Berkenkamp et~al.(2017)Berkenkamp, Turchetta, Schoellig, and
  Krause]{berkenkamp2017safe}
F.~Berkenkamp, M.~Turchetta, A.~Schoellig, and A.~Krause.
\newblock Safe model-based reinforcement learning with stability guarantees.
\newblock In \emph{Advances in {N}eural {I}nformation {P}rocessing {S}ystems},
  pages 908--918, 2017.

\bibitem[Bertsekas(2005)]{bertsekas2005dynamic}
D.~P. Bertsekas.
\newblock Dynamic programming and suboptimal control: A survey from {ADP} to
  {MPC}.
\newblock \emph{European Journal of Control}, 11\penalty0 (4-5):\penalty0
  310--334, 2005.

\bibitem[Bertsekas(2012)]{Bertsekas:12}
D.~P. Bertsekas.
\newblock \emph{Dynamic Programming and Optimal Control}, volume~2.
\newblock Athena Scientific, Nashua, USA, 4th edition, 2012.

\bibitem[{Bertsekas}(2017)]{Bertsekas:TNNLS}
D.~P. {Bertsekas}.
\newblock Value and policy iterations in optimal control and adaptive dynamic
  programming.
\newblock \emph{IEEE Transactions on Neural Networks and Learning Systems},
  28\penalty0 (3):\penalty0 500--509, 2017.
\newblock \doi{10.1109/TNNLS.2015.2503980}.

\bibitem[Bian and Jiang(2016)]{BIAN2016348}
T.~Bian and Z.-P. Jiang.
\newblock Value iteration and adaptive dynamic programming for data-driven
  adaptive optimal control design.
\newblock \emph{Automatica}, 71:\penalty0 348 -- 360, 2016.
\newblock ISSN 0005-1098.
\newblock \doi{https://doi.org/10.1016/j.automatica.2016.05.003}.

\bibitem[Buşoniu et~al.(2018)Buşoniu, {de Bruin}, Tolić, Kober, and
  Palunko]{BUSONIU20188}
L.~Buşoniu, T.~{de Bruin}, D.~Tolić, J.~Kober, and I.~Palunko.
\newblock Reinforcement learning for control: Performance, stability, and deep
  approximators.
\newblock \emph{Annual Reviews in Control}, 46:\penalty0 8 -- 28, 2018.
\newblock ISSN 1367-5788.
\newblock \doi{https://doi.org/10.1016/j.arcontrol.2018.09.005}.

\bibitem[Granzotto(2019)]{granz2019}
M.~Granzotto.
\newblock \emph{Near-optimal control of discrete-time nonlinear systems with
  stability guarantees}.
\newblock PhD thesis, Université de Lorraine, 2019.
\newblock URL \url{http://www.theses.fr/2019LORR0301}.

\bibitem[Granzotto et~al.(2019)Granzotto, Postoyan, Buşoniu, Nešić, and
  Daafouz]{granzottoCDC2019}
M.~Granzotto, R.~Postoyan, L.~Buşoniu, D.~Nešić, and J.~Daafouz.
\newblock Optimistic planning for the near-optimal control of nonlinear
  switched discrete-time systems with stability guarantees.
\newblock In \emph{IEEE Conference on Decision and Control}, Nice, France,
  2019.
\newblock URL \url{https://arxiv.org/pdf/1908.01404.pdf}.

\bibitem[Granzotto et~al.(2020{\natexlab{a}})Granzotto, Postoyan, Buşoniu,
  Nešić, and Daafouz]{granzotto2020}
M.~Granzotto, R.~Postoyan, L.~Buşoniu, D.~Nešić, and J.~Daafouz.
\newblock Finite-horizon discounted optimal control: stability and performance.
\newblock \emph{IEEE Transactions on Automatic Control}, 2020{\natexlab{a}}.
\newblock \doi{10.1109/TAC2020.2985904}.

\bibitem[Granzotto et~al.(2020{\natexlab{b}})Granzotto, Postoyan, Buşoniu,
  Nešić, and Daafouz]{granzotto2020optimistic}
M.~Granzotto, R.~Postoyan, L.~Buşoniu, D.~Nešić, and J.~Daafouz.
\newblock Stable near-optimal control of nonlinear switched discrete-time
  systems: a planning-based approach.
\newblock In \emph{Submitted to journal publication}, 2020{\natexlab{b}}.

\bibitem[Grimm et~al.(2005)Grimm, Messina, Tuna, and Teel]{grimm2005}
G.~Grimm, M.~J. Messina, S.~E. Tuna, and A.~R. Teel.
\newblock Model predictive control: for want of a local control {Lyapunov}
  function, all is not lost.
\newblock \emph{IEEE Transactions on Automatic Control}, 50\penalty0
  (5):\penalty0 546--558, 2005.
\newblock ISSN 0018-9286.
\newblock \doi{10.1109/TAC.2005.847055}.

\bibitem[{Gr\"une} and Rantzer(2008)]{gruneperformance}
L.~{Gr\"une} and A.~Rantzer.
\newblock On the infinite horizon performance of receding horizon controllers.
\newblock \emph{IEEE Transactions on Automatic Control}, 53\penalty0
  (9):\penalty0 2100--2111, 2008.
\newblock ISSN 0018-9286.
\newblock \doi{10.1109/TAC.2008.927799}.

\bibitem[{Heydari}(2014)]{heydari2014adp}
A.~{Heydari}.
\newblock Revisiting approximate dynamic programming and its convergence.
\newblock \emph{IEEE Transactions on Cybernetics}, 44\penalty0 (12):\penalty0
  2733--2743, 2014.
\newblock \doi{10.1109/TCYB.2014.2314612}.

\bibitem[{Heydari}(2016)]{heydari2016acc}
A.~{Heydari}.
\newblock Analysis of stabilizing value iteration for adaptive optimal control.
\newblock In \emph{2016 American Control Conference (ACC)}, pages 5746--5751,
  2016.
\newblock \doi{10.1109/ACC.2016.7526570}.

\bibitem[Heydari(2017)]{heydari2017stability}
A.~Heydari.
\newblock Stability analysis of optimal adaptive control under value iteration
  using a stabilizing initial policy.
\newblock \emph{IEEE Transactions on Neural Networks and Learning Systems},
  29\penalty0 (9):\penalty0 4522--4527, 2017.

\bibitem[Heydari(2018)]{heydari}
A.~Heydari.
\newblock Stability analysis of optimal adaptive control using value iteration
  with approximation errors.
\newblock \emph{IEEE Transactions on Automatic Control}, 2018.
\newblock ISSN 0018-9286.
\newblock \doi{10.1109/TAC.2018.2790260}.

\bibitem[Hren and Munos(2008)]{hren2008}
J.-F. Hren and R.~Munos.
\newblock Optimistic planning of deterministic systems.
\newblock In \emph{European Workshop on Reinforcement Learning}, pages
  151--164, Villeneuve d'Ascq, France, 2008.

\bibitem[Jiang and Jiang(2012)]{JIANG20122699}
Y.~Jiang and Z.-P. Jiang.
\newblock Computational adaptive optimal control for continuous-time linear
  systems with completely unknown dynamics.
\newblock \emph{Automatica}, 48\penalty0 (10):\penalty0 2699 -- 2704, 2012.
\newblock ISSN 0005-1098.
\newblock \doi{https://doi.org/10.1016/j.automatica.2012.06.096}.

\bibitem[Keerthi and Gilbert(1985)]{keerthi1985}
S.~Keerthi and E.~Gilbert.
\newblock An existence theorem for discrete-time infinite-horizon optimal
  control problems.
\newblock \emph{IEEE Transactions on Automatic Control}, 30\penalty0
  (9):\penalty0 907--909, 1985.
\newblock ISSN 0018-9286.
\newblock \doi{10.1109/TAC.1985.1104084}.

\bibitem[Kiumarsi et~al.(2017)Kiumarsi, Lewis, and Jiang]{kiumarsi2017h}
B.~Kiumarsi, F.~L. Lewis, and Z.-P. Jiang.
\newblock {H$\infty$} control of linear discrete-time systems: Off-policy
  reinforcement learning.
\newblock \emph{Automatica}, 78:\penalty0 144--152, 2017.

\bibitem[{Lewis} and {Vrabie}(2009)]{LewisRLDP2009}
F.~L. {Lewis} and D.~{Vrabie}.
\newblock Reinforcement learning and adaptive dynamic programming for feedback
  control.
\newblock \emph{IEEE Circuits and Systems Magazine}, 9\penalty0 (3):\penalty0
  32--50, 2009.
\newblock \doi{10.1109/MCAS.2009.933854}.

\bibitem[{Liu} et~al.(2015){Liu}, {Li}, and {Wang}]{derongliu2015}
D.~{Liu}, H.~{Li}, and D.~{Wang}.
\newblock Error bounds of adaptive dynamic programming algorithms for solving
  undiscounted optimal control problems.
\newblock \emph{IEEE Transactions on Neural Networks and Learning Systems},
  26\penalty0 (6):\penalty0 1323--1334, 2015.
\newblock \doi{10.1109/TNNLS.2015.2402203}.

\bibitem[Pang et~al.(2019)Pang, Bian, and Jiang]{pang2019adaptive}
B.~Pang, T.~Bian, and Z.-P. Jiang.
\newblock Adaptive dynamic programming for finite-horizon optimal control of
  linear time-varying discrete-time systems.
\newblock \emph{Control Theory and Technology}, 17\penalty0 (1):\penalty0
  73--84, 2019.

\bibitem[{Pavlov} et~al.(2019){Pavlov}, {Shames}, and {Manzie}]{Pavlov2019}
A.~{Pavlov}, I.~{Shames}, and C.~{Manzie}.
\newblock Early termination of {NMPC} interior point solvers: Relating the
  duality gap to stability.
\newblock In \emph{2019 18th European Control Conference (ECC)}, pages
  805--810, 2019.
\newblock \doi{10.23919/ECC.2019.8795629}.

\bibitem[Postoyan et~al.(2017)Postoyan, Buşoniu, Nešić, and
  Daafouz]{romain2016}
R.~Postoyan, L.~Buşoniu, D.~Nešić, and J.~Daafouz.
\newblock Stability analysis of discrete-time infinite-horizon optimal control
  with discounted cost.
\newblock \emph{IEEE Transactions on Automatic Control}, 62\penalty0
  (6):\penalty0 2736--2749, 2017.
\newblock ISSN 0018-9286.
\newblock \doi{10.1109/TAC.2016.2616644}.

\bibitem[Postoyan et~al.(2019)Postoyan, Granzotto, Buşoniu, Scherrer, Nešić,
  and Daafouz]{romainAVICDC2019}
R.~Postoyan, M.~Granzotto, L.~Buşoniu, B.~Scherrer, D.~Nešić, and
  J.~Daafouz.
\newblock Stability guarantees for nonlinear discrete-time systems controlled
  by approximate value iteration.
\newblock In \emph{IEEE Conference on Decision and Control}, Nice, France,
  2019.

\bibitem[Sutton and Barto(2017)]{Sutton}
R.~S. Sutton and A.~G. Barto.
\newblock \emph{Reinforcement Learning: An Introduction}.
\newblock MIT Press, Cambridge, USA, 2nd edition, 2017.

\bibitem[Wei et~al.(2015)Wei, Liu, and Lin]{wei2015value}
Q.~Wei, D.~Liu, and H.~Lin.
\newblock Value iteration adaptive dynamic programming for optimal control of
  discrete-time nonlinear systems.
\newblock \emph{IEEE Transactions on Cybernetics}, 46\penalty0 (3):\penalty0
  840--853, 2015.

\end{thebibliography}

\end{document}